\documentclass[12pt]{amsart}  
\usepackage{amsthm, amscd}
\usepackage{amssymb,latexsym}
\input xy
\xyoption{all} 

\theoremstyle{plain}
\newtheorem*{introtheorem}{Theorem}
\newtheorem*{introcorollary}{Corollary}

\newtheorem{theorem}{Theorem}[section]
\newtheorem{proposition}[theorem]{Proposition}
\newtheorem{corollary}[theorem]{Corollary}
\newtheorem{lemma}[theorem]{Lemma}
\newtheorem{conjecture}[theorem]{Conjecture}
 
\theoremstyle{definition}    
\newtheorem{definition}[theorem]{Definition}
\newtheorem{remark}[theorem]{Remark}
\newtheorem{example}[theorem]{Example}

\theoremstyle{remark}

\newcommand{\ZZ}{\mathbf{Z}}

\newcommand{\OO}{\mathcal{O}}
\newcommand{\II}{\mathcal{I}}
   
\newcommand{\FF}{\mathcal{F}}
\newcommand{\GG}{\mathcal{G}}
\newcommand{\LL}{\mathcal{L}}

\begin{document}
 
\title{\bf Verlinde bundles and generalized theta linear series}
\author{Mihnea Popa}
\address{Department of Mathematics
\\University of Michigan \\ 525 East University \\ Ann Arbor, MI 48109-1109  and Institute of
Mathematics of the Romanian Academy}
\email{mpopa@math.lsa.umich.edu}

\subjclass{Primary 14H60; Secondary 14J60}
\keywords{Vector bundles, nonabelian theta functions}
\maketitle

$$\textbf{Introduction}$$

The purpose of this paper is to introduce and study some very interesting vector bundles 
on the Jacobian of a curve which are associated to
generalized theta linear series on moduli spaces of vector 
bundles on that curve. 
By way of application, we establish effective bounds for the global and normal generation 
of such linear series and we give new proofs of the known results and a new perspective 
on the familiar conjectures concerning duality for generalized theta functions. We also 
show that these bundles lead to new examples (in the spirit of \cite{Raynaud}) of base 
points for the determinant linear series.

Let $X$ be a smooth projective complex curve of genus $g\geq 2$. 
Let also $U_{X}(r,0)$ be the moduli space of semistable vector bundles of rank $r$ 
and degree $0$ on $X$ and $SU_{X}(r)$ the moduli space of semistable rank $r$ vector 
bundles with trivial determinant.
We define the \emph{Verlinde bundles} to be push-forwards of  
pluritheta line bundles on $U_{X}(r,0)$ to the Jacobian of $X$ by the determinant 
map ${\rm det}:U_{X}(r,0)\rightarrow J(X)$:
$$E_{r,k}:={\rm det}_{*}\OO(k\Theta_{N}),$$
where $\Theta_{N}$ is the generalized theta divisor associated to a line bundle 
$N\in {\rm Pic}^{g-1}(X)$ (see \emph{1.1}). Their fibers are precisely the Verlinde 
vector spaces of level $k$ theta functions $H^{0}(SU_{X}(r),\LL^{k})$, where $\LL$ is 
the determinant line bundle on $SU_{X}(r)$.
These bundles are the focus of the present paper. Our viewpoint is that they allow a translation
of statements about spaces of generalized theta functions into global statements about 
vector bundles on Jacobians. This in turn opens the way towards a more geometric study of 
these spaces.

We study at some length the behavior of the Verlinde bundles under natural operations associated to 
vector bundles on abelian varieties, 
e.g. in the spirit 
of Mukai \cite{Mukai}.
In particular, we focus on their 
Fourier-Mukai transforms $\widehat{E_{r,k}}$.
An interesting picture is provided by the relationship between the bundles $E_{r,k}$ and 
$E_{k,r}$ via the Fourier transform. Identifying $J(X)$ with the dual ${\rm Pic}^{0}(J(X))$ 
in the canonical way, we prove:

\begin{introtheorem}
There is a (non-canonical) isomorphism $r_{J}^{*}E_{r,k}^{*}\cong r_{J}^{*}\widehat{E_{k,r}}$, where $r_{J}$ 
is the multiplication by $r$ map on $J(X)$. In particular we have (cf. \cite{Donagi}):
$$h^{0}(SU_{X}(k),\LL^{r})\cdot r^{g}=h^{0}(U_{X}(k,0),\OO(r\Theta_{N}))
\cdot k^{g}.$$
\end{introtheorem}

\noindent
The last statement is the duality result of \cite{Donagi} on dimensions of spaces of generalized 
theta functions. It is natural then, as in the case of theta functions, to try to understand stronger  
dualities between these bundles. The main result in this direction is:

\begin{introtheorem}  
For every $k\geq 2$, there is a canonical nontrivial map
$$SD:E_{r,k}^{*}\longrightarrow \widehat{E_{k,r}}.$$
In the case $k=1$, $E_{r,1}$ and $\widehat
{E_{1,r}}$ are stable and $SD$ is an isomorphism.
\end{introtheorem}

\noindent
The last part follows once we prove a statement of independent interest concerning the stability 
of the Fourier-Mukai transform of an arbitrary polarization on an abelian variety.
Note also that we are granting the Verlinde formula in the proof of the theorem.

We recall that in the context of nonabelian theta functions, one of the central problems is a conjectural 
duality between certain spaces of $GL(n)$ and $SL(n)$-theta functions, known as the 
\emph{strange duality conjecture}. It asserts that there is a canonical isomorphism
$$H^{0}(SU_{X}(r),\LL^{k})^{*}\overset{\cong}{\longrightarrow} H^{0}(U_{X}(k,0),\OO(r\Theta_{N})).$$
These two vector spaces can be seen as fibers of the vector bundles $E_{r,k}^{*}$ and $\widehat{E_{k,r}}$
respectively, and indeed the natural map $SD$ in the theorem induces the strange duality map on 
the fibers. The conjecture is consequently equivalent to the global statement saying that 
$SD$ is an isomorphism for every $k$. 
In particular, the theorem already implies the
\emph{strange duality at level $1$}, initially proved in \cite{BNR}, saying that
$$H^{0}(SU_{X}(r),\LL)^{*}\overset{SD}{\longrightarrow}
H^{0}(J(X),\OO(r\Theta_{N}))$$ 
is an isomorphism. Moreover, the first theorem above can also be seen as providing some global 
evidence for the general conjecture.
One expects that an understanding of further properties 
characterizing the Verlinde bundles and the 
maps between them (e.g. as in Remark \ref{51}) could shed some new light on this whole
circle of ideas.

Our main concrete application of the Verlinde bundles is to the  
study of effective global generation and normal generation for pluritheta 
line bundles on $U_{X}(r,0)$. 
The underlying idea is very simple: the determinant map makes 
$U_{X}(r,0)$ a fiber space over $J(X)$ and one expects that effective bounds 
on the base and the fibers should give bounds on the total space.
For the fiberwise situation, i.e. the case of multiples of the determinant bundle on 
$SU_{X}(r)$, we make use of recent results proved in \cite{Popa2}. Thus the complications 
arise when one tries to understand the global or normal generation of $E_{r,k}$ on $J(X)$.
We begin by proving that $E_{r,k}$ is globally generated 
if and only if  $k\geq r+1$.
The consequence of this fact for linear series is stated below
in a form that allows algorithmic applications. 

\begin{introtheorem}
$\OO(k\Theta_{N})$ is globally generated on $U_{X}(r,0)$ as long as $k\geq r+1$
and $\LL^{k}$ is globally generated on $SU_{X}(r)$. Moreover, $\OO(r\Theta_{N})$ is 
not globally generated.
\end{introtheorem}
\noindent
As mentioned above, this provides effective results when combined with the bound obtained 
in \cite{Popa2} \S4 for $SU_{X}(r)$.

\begin{introcorollary}
$\OO(k\Theta_{N})$ is globally generated on $U_{X}(r,0)$ if 
$k\geq \frac{(r+1)^{2}}{4}$.
\end{introcorollary}
\noindent
Note that the theorem suggests that the optimal result on $U_{X}(r,0)$ could be the global generation 
of $\OO((r+1)\Theta_{N})$, which goes beyond what we currently have on $SU_{X}(r)$. This is 
indeed true for moduli of rank $2$ and rank $3$ vector bundles (see \ref{60}) and accounts for 
the discrepancy between the way we are formulating the theorem and the statement of the corollary.  
For some general results and questions in the direction of optimal effective base point 
freeness the reader can consult \cite{Popa2} \S4 and \S5. 

Moving to effective surjectivity of multiplication maps, we prove the following:
\begin{introtheorem}
(a) The multiplication map 
$$H^{0}(E_{r,k})\otimes H^{0}(E_{r,k})\longrightarrow H^{0}(E_{r,k}^{\otimes 2})$$
is surjective for $k\geq 2r+1$.
\newline
(b) Assume that $E_{r,k}$ on $J(X)$ and $\LL^{k}$ on $SU_{X}(r)$ are globally generated. 
Then the multiplication map 
$$\mu_{k}:H^{0}(\OO(k\Theta_{N}))\otimes 
H^{0}(\OO(k\Theta_{N}))\longrightarrow H^{0}(\OO(2k\Theta_{N}))$$ 
is surjective as long as the multiplication map $H^{0}(\LL^{k})\otimes H^{0}(\LL^{k})
\rightarrow H^{0}(\LL^{2k})$ on $SU_{X}(r)$ is surjective and $k\geq 2r+1$. 
In particular, for such $k$, $\OO(k\Theta_{N})$ is very ample.
\end{introtheorem} 
\noindent
To get an effective application in the case of rank $2$ vector bundles one can use a theorem
of Laszlo \cite{Laszlo}. It is interesting to note that while the very ampleness of the 
determinant line bundle on $SU_{X}(2)$ has been the subject of various approaches (see e.g. 
\cite{Verra}, \cite{Izadi}, \cite{Laszlo}), no effective results seem to have been known for generalized theta
bundles on $U_{X}(2,0)$.

\begin{introcorollary}
For a generic curve $X$, the multiplication map 
$$H^{0}(\OO(k\Theta_{N}))\otimes H^{0}(\OO(k\Theta_{N}))\overset{\mu_{k}}
{\longrightarrow}H^{0}(\OO(2k\Theta_{N}))$$ 
on $U_{X}(2,0)$ is surjective for $k\geq {\rm max}\{5,g-2\}$ and so $\OO(k\Theta_{N})$ 
is very ample for such $k$. 
\end{introcorollary}

\noindent
The proofs of both theorems make use of some very interesting
recent results of Pareschi \cite{Pareschi} (see \emph{1.4}), drawing also on 
earlier work of Kempf, on cohomological criteria 
for global generation and normal generation of vector bundles on abelian varieties 
$\footnote{Pareschi uses these criteria in \cite{Pareschi} to prove a conjecture of Lazarsfeld 
on syzygies of abelian varieties.}$.

A somewhat surprising
application of the Verlinde bundles is to the problem of finding base points 
for the linear system $|\LL|$ on $SU_{X}(r)$ (see \cite{Beauville} \S2 and \cite{Popa} 
for a survey and results in this direction).
We extend a method
of Raynaud \cite{Raynaud} to show that the Fourier-Mukai transforms 
$\widehat{E_{r,k}^{*}}$, restricted to certain 
embeddings of $X$ in $J(X)$, provide new examples of base points for  
$|\LL|$. Their ranks are given by the Verlinde formula at all levels $k\geq 1$.  

\begin{introtheorem}
There exists an embedding $X\hookrightarrow J(X)$ such that $\widehat{E_{r,k}^{*}}_{|X}$
gives (via twisting by a line bundle) a base point for the linear system $|\LL|$ 
on $SU_{X}(s_{k,r})$, where $s_{k,r}:=h^{0}(SU_{X}(k), \LL^{r})$. 
\end{introtheorem}
\noindent
This construction turns out to be a natural generalization of the examples given 
in \cite{Raynaud}, which in our language correspond to the level $1$ Verlinde bundles.

Most of our results can be formulated in the setting of moduli of vector
bundles of arbitrary rank and degree. To keep the arguments uniform and at a technical 
minimum we have chosen to explain all the details for the moduli spaces $U_{X}(r,0)$ and
$SU_{X}(r)$.
However, the last section contains complete statements of the analogous results that hold 
in the general case, with hints for the (minor) modifications needed in the proofs. 

We briefly explain the structure of the paper, which is mainly ordered according to the increasing 
complexity of the methods involved. The first section provides the background material needed, with rather detailed 
statements included for the convenience of the reader. The definition and basic study of the Verlinde 
bundles occupy \S2. The third section is an account of the Fourier duality picture for Verlinde bundles 
with applications to duality for generalized theta functions. In \S4 we explain the construction of 
base points for the determinant linear series. The fifth section contains the main application to the effective 
global and normal generation problem for pluritheta linear series. Finally, in the last section we 
explain how the results can be extended to the setting of arbitrary degree vector bundles. 

There are a few very interesting related topics that have been left out of this paper.
The reader familiar with the theory of nonabelian theta functions will notice that there is 
a certain parallel between the particular form of the Verlinde bundles $E_{r,k}$ and 
(in a vague sense) the representation
theory of the Verlinde vector spaces $H^{0}(SU_{X}(r),\LL^{k})$. In general one studies the   
representation theory of the Verlinde spaces in order to learn something about 
linear series and the geometry of the moduli 
spaces. Our hope is that a better understanding of the global properties of the $E_{r,k}$'s
might, in the other direction, say something about the representations.
From a different point of view, the theory developed here may be at least partially generalized 
to the setting of moduli of sheaves on smooth (irregular) surfaces. We are hoping  to address these 
questions in future work. 

\textbf{Acknowledgements.} I would like to thank my advisor, R. Lazarsfeld, for valuable 
suggestions and 
support throughout this project. My thanks also go to I. Dolgachev for numerous fruitful 
conversations on the subject and to A. Beauville for discussions at the 1999 VBAC conference. 
Finally I am indebted to G. Pareschi for allowing my access to 
various versions of his work \cite{Pareschi}; this has been essential for completing the 
present paper.

\section{\textbf{Preliminaries and notations}}

\emph{1.1. Generalized theta divisors.} 
Let $U_{X}(r,0)$ be the moduli space of equivalence classes of semistable vector 
bundles of rank $r$ and degree $0$ on $X$, $SU_{X}(r,L)$ the moduli space of vector 
bundles with fixed determinant $L\in$ Pic$^{0}(X)$ and $U_{X}^{s}(r,0)$ and 
$SU_{X}^{s}(r,L)$ the open subsets corresponding to stable bundles. We recall the 
construction and some basic facts       
about generalized theta divisors on these 
moduli spaces, drawing especially on \cite{Drezet}. Analogous constructions
work for any degree $d$, as we will recall in the last section.

Let $N$ be a line bundle of degree $g-1$ on $X$. Then $\chi(E\otimes N)=0$ for all 
$E\in U_{X}(r,0)$. Consider the subset of $U_{X}^{s}(r,0)$
$$\Theta_{N}^{s}:= \{E\in U_{X}^{s}(r,0)|~ h^{0}(E\otimes N)\neq 0\}$$
and the analogous set in $SU_{X}^{s}(r,L)$. One can prove (see \cite{Drezet} (7.4.2))
that $\Theta_{N}^{s}$ is a hypersurface in $U_{X}^{s}(r,0)$ (resp. $SU_{X}^{s}(r,L)$).
Denote by $\Theta_{N}$ the closure of $\Theta_{N}^{s}$ in $U_{X}(r,0)$ and 
$SU_{X}(r,L)$. As we vary $N$, these hypersurfaces are called $\it{generalized~ theta~
divisors}$. It is proved in \cite{Drezet}, Theorem A, that $U_{X}(r,0)$ and $SU_{X}(r,L)$
are locally factorial and so the generalized theta divisors determine line bundles
$\OO(\Theta_{N})$ on these moduli spaces. We have the following important facts:

\begin{theorem}(\cite{Drezet}, Theorem B)  
The line bundle $\OO(\Theta_{N})$ on $SU_{X}(r,L)$ does not depend on the choice
of $N$. The Picard group of $SU_{X}(r,L)$ is isomorphic to $\ZZ$, generated by 
$\OO(\Theta_{N})$. 
\end{theorem}
\medskip
The line bundle in the theorem above, independent of the choice of $N$, is denoted 
by $\LL$ and is called the $\it{determinant}$ $\it{bundle}$.
\newline
\begin{theorem}(\cite{Drezet}, Theorem C)\label{47}
The inclusions ${\rm Pic}(J(X))\subset {\rm Pic}(U_{X}(r,0))$ (given by the determinant morphism)
and $\ZZ\cdot\OO(\Theta_{N})\subset {\rm Pic}(U_{X}(r,0))$ induce an isomorphism
$${\rm Pic}(U_{X}(r,0))\cong {\rm Pic}(J(X))\oplus\ZZ.$$
\end{theorem} 

More generally, for any vector bundle $F$ of rank $k$ and degree $k(g-1)$, we can 
define $$\Theta_{F}^{s}:= \{E\in U_{X}^{s}(r,0)|~ h^{0}(E\otimes F)\neq 0\}$$
and denote by $\Theta_{F}$ the closure of $\Theta_{F}^{s}$ in $U_{X}(r,0)$. 
It is clear that, for generic $F$ at least, $\Theta_{F}$ is 
strictly contained in $U_{X}(r,0)$ (in which case it is again a divisor). It is useful
to know what is the dependence of $\OO(\Theta_{F})$ on $F$ and in this direction 
we have:

\begin{proposition}(\cite{Drezet} (7.4.3) and \cite{Donagi} Prop.3)\label{32}
Let $F$ and $G$ be two vector bundles of slope $g-1$ on $X$. If ${\rm rk}(F)= m\cdot 
{\rm rk}(G)$, then 
$$\OO(\Theta_{F})\cong \OO(\Theta_{G})^{\otimes m}\otimes {\rm det}^{*}({\rm det}F\otimes 
({\rm detG})^{\otimes -m})$$
where we use the natural identification of ${\rm Pic}^{0}(X)$ with ${\rm Pic}^{0}(J(X))$.
In particular, if $F$ has rank $k$ and $N$ is a line bundle of degree $g-1$,
we get
$$\OO(\Theta_{F})\cong \OO(\Theta_{N})^{\otimes k}\otimes {\rm det}^{*}({\rm det}F\otimes 
N^{\otimes -k}).$$
\end{proposition}

\emph{1.2. A convention on theta divisors.}
When looking at generalized theta divisors $\Theta_{N}$ and the corresponding linear 
series, it will be convenient to consider the line bundle $N$ to be a $\it{theta ~
characteristic}$, i.e. satisfying $N^{\otimes 2}\cong \omega_{X}$. The assumption 
brings some simplifications to most of the arguments, but on the other hand this 
case implies all the results for arbitrary $N$. This is true since for any $N$ and $M$
in Pic$^{g-1}(X)$, if $\xi:=N\otimes M^{-1}$, twisting by $\xi$ gives an automorphism:
$$U_{X}(r,0)\overset{\otimes\xi}{\longrightarrow}U_{X}(r,0)$$
by which $\OO(\Theta_{M})$ corresponds to $\OO(\Theta_{N})$. 
As an example, we will use freely isomorphisms of the type $r_{J}^{*}\OO_{J}(\Theta_{N})
\cong \OO_{J}(r^{2}\Theta_{N})$, where $r_{J}$ is multiplication by $r$ on $J(X)$, which
in general would work only up to numerical equivalence.  
One can easily see that in each 
particular proof the arguments could be worked out in the general situation with little 
extra effort (the main point is that the cohomological arguments work even if we use numerical
equivalence instead of linear equivalence). 

It is worth mentioning another convention about the notation that we will be using. The 
divisors $\Theta_{N}$ make sense of course on both $U_{X}(r,0)$ for $r\geq 2$ and 
$J(X)=U_{X}(1,0)$ and in some proofs both versions will be used. We will denote the 
associated line bundle simply by $\OO(\Theta_{N})$ if $\Theta_{N}$ lives on $U_{X}(r,0)$
and by $\OO_{J}(\Theta_{N})$ if it lives on the Jacobian. 

\emph{1.3. The Fourier-Mukai transform on an abelian variety.} 
Here we give a brief overview of some basic facts on the Fourier-Mukai transform on 
an abelian variety, following the original paper of Mukai \cite{Mukai}.
Let $X$ be an abelian variety of dimension $g$, $\widehat{X}$ its dual and 
$\mathcal{P}$ the Poincar\'e
line bundle on $X\times \widehat{X}$, normalized such that 
$\mathcal{P}_{|X\times\{0\}}$
and $\mathcal{P}_{|\{0\}\times \widehat{X}}$ are trivial. To any coherent sheaf $\FF$ 
on $X$ we can associate the sheaf ${p_{2}}_{*}({p_{1}}^{*}\FF\otimes\mathcal{P})$ on
$\widehat{X}$ via the natural diagram:
$$\xymatrix{ & & X\times \widehat{X} \ar[dl]^{p_{1}} \ar[dr]^{p_{2}} \\
& X & & \widehat{X} } $$
This correspondence gives a functor
$$\mathcal{S}: {\rm Coh(X)}\rightarrow {\rm Coh}(\widehat{X}).$$ 
If we denote by $D(X)$ and 
$D(\widehat{X})$ the derived categories of Coh$(X)$ and Coh$(\widehat{X})$, then 
the derived functor $\mathbf{R}\mathcal{S}:D(X)\rightarrow D(\widehat{X})$ is defined 
(and called the Fourier functor)
and one can consider $\mathbf{R}\widehat{\mathcal{S}}:D(\widehat{X})\rightarrow D(X)$ 
in a similar way. Mukai's main theorem is the following:

\begin{theorem}(\cite{Mukai} (2.2))\label{34}
The Fourier functor establishes an equivalence of categories between $D(X)$ and 
$D(\widehat{X})$. More precisely there are isomorphisms of functors:
$$\mathbf{R}\mathcal{S}\circ\mathbf{R}\widehat{\mathcal{S}}\cong (-1_{\widehat{X}})^{*}
[-g]$$
$$\mathbf{R}\widehat{\mathcal{S}}\circ\mathbf{R}\mathcal{S}\cong (-1_{X})^{*}
[-g].$$
\end{theorem}
In this paper we will essentially have to deal with the simple situation when by 
applying the Fourier functor we get back another vector bundle, i.e. a complex with 
only one nonzero (and locally free) term. This is packaged in the following 
definition (see \cite{Mukai} (2.3):
\begin{definition}
A coherent sheaf $\FF$ on $X$ satisfies I.T. (index theorem) with $\it{index ~j}$ if 
$H^{i}(\FF\otimes\alpha)=0$ for all $\alpha\in {\rm Pic}^{0}(X)$ and all $i\neq j$.
In this situation we have $R^{i}\mathcal{S}(\FF)=0$ for all $i\neq j$ 
and by the base change theorem $R^{j}\mathcal{S}(\FF)$ is locally free. We 
denote $R^{j}\mathcal{S}(\FF)$ by $\widehat{\FF}$ and call it the $\it{Fourier~
transform}$ of $\FF$. Note that then $\mathbf{R}\mathcal{S}(\FF)\cong \widehat{F}[-j]$.
\end{definition}
For later reference we list some basic properties of the Fourier transform that will 
be used repeatedly throughout the paper:

\begin{proposition}\label{30}
Let $X$ be an abelian variety and $\widehat{X}$ its dual and let $\mathbf{R}\mathcal
{S}$ and $\mathbf{R}\widehat{\mathcal{S}}$ the corresponding Fourier functors. Then 
the following are true:

(1)(\cite{Mukai} (3.4)) Let $Y$ be an abelian veriety, $f:Y\rightarrow X$ an isogeny 
and $\hat{f}:\widehat{X}\rightarrow \widehat{Y}$ the dual isogeny. Then there are 
isomorphisms of functors:
$$f^{*}\circ\mathbf{R}\widehat{\mathcal{S}}_{X}\cong \mathbf{R}\widehat{\mathcal{S}}
_{Y}\circ \hat{f}_{*}$$
$$f_{*}\circ\mathbf{R}\widehat{\mathcal{S}}_{Y}\cong \mathbf{R}\widehat{\mathcal{S}}
_{X}\circ \hat{f}^{*}.$$

(2)(\cite{Mukai} (3.7)) Let $\FF$ and $\GG$ coherent sheaves on $X$ and define their 
Pontrjagin product by $\FF *\GG := \mu_{*}({p_{1}}^{*}\FF\otimes {p_{2}}^{*}\GG)$,
where $\mu: X\times X\rightarrow X$ is the multiplication on $X$. Then we have the 
following isomorphisms:
$$\mathbf{R}\mathcal{S}(\FF *\GG)\cong \mathbf{R}\mathcal{S}(\FF)\otimes \mathbf{R}
\mathcal{S}(\GG)$$
$$\mathbf{R}\mathcal{S}(\FF\otimes\GG)\cong \mathbf{R}\mathcal{S}(\FF) *\mathbf{R}
\mathcal{S}(\GG)[g],$$ 
where the operations on the right hand side should be thought of in the derived 
category.

(3)(\cite{Mukai} (3.11)) Let $L$ be a nondegenerate line bundle on $X$ of index $i$,
i.e. $h^{i}(L)\neq 0$ and $h^{j}(L)=0$ for all $j\neq i$.
Then by \cite{Mumford2} \S16, I.T. holds for $L$ and there is an isomorphism
$$\phi_{L}^{*}\widehat{L}\cong H^{i}(L)\otimes L^{-1}~(\cong \underset{|\chi (L)|}
{\bigoplus}L^{-1}),$$ with $\phi_{L}$
the isogeny canonically defined by {L}.
 
(4)(\cite{Mukai} (3.1)) Let $\FF$ be a coherent sheaf on $X$, $x\in \widehat{X}$
and $P_{x}\in {\rm Pic}^{0}(X)$ the corresponding line bundle. Then we have an 
isomorphism:
$$\mathbf{R}\mathcal{S}(\FF\otimes P_{x})\cong t_{x}^{*}\mathbf{R}\mathcal{S}(\FF),$$
where $t_{x}$ is translation by $x$. 

(5)(\cite{Mukai} (2.4),(2.5),(2.8)) Assume that $E$ satisfies I.T. with index $i$.
Then $\widehat{E}$ satisfies I.T. with index $g-i$ and in this case 
$$\chi(E)=(-1)^{i}\cdot {\rm rk}(\widehat{E}).$$
Moreover, there are isomorphisms ${\rm Ext}^{k}(E,E)\cong 
{\rm Ext}^{k}(\widehat{E},\widehat{E})$
for all $k$, and in particular $E$ is simple if and only if $\widehat{E}$ is simple.
\end{proposition}

\emph{1.4. Global generation and normal generation of vector bundles on abelian 
varieties.}
For the reader's convenience, in the present paragraph we give a brief account of 
some very recent results and techniques in the study of vector bundles on abelian 
varieties -- following work of Pareschi \cite{Pareschi} -- in a form convenient for
our purposes. The underlying theme
is to give useful criteria for the global generation and surjectivity of
multiplication maps of such vector bundles.

Let $X$ be an abelian variety and $E$ a vector bundle on $X$. Building on earlier 
work of Kempf \cite{Kempf}, Pareschi proves the following cohomological criterion 
for global generation:

\begin{theorem}(\cite{Pareschi} (2.1))\label{40}
Assume that $E$ satisfies the following vanishing property:
$$h^{i}(E\otimes\alpha)=0, ~ \forall\alpha\in{\rm Pic}^{0}(X)~{\rm and} ~ \forall 
i>0.$$
Then for any ample line bundle $L$ on $X$, $E\otimes L$ is globally generated.
\end{theorem}

In another direction, in order to attack questions about multiplication maps of the
form
\begin{equation}
H^{0}(E)\otimes H^{0}(F)\longrightarrow H^{0}(E\otimes F)
\end{equation}
for $E$ and $F$ vector bundles on $X$, the right notion turns out to be that of 
skew Pontrjagin product:

\begin{definition}\label{80}
Let $E$ and $F$ be coherent sheaves on the abelian variety $X$. Then the 
$\it{skew ~Pontrjagin ~product}$ of $E$ and $F$ is defined by 
$$E\hat{*} F := {p_{1}}_{*}((p_{1}+p_{2})^{*}E\otimes p_{2}^{*}F)$$
where $p_{1},p_{2}:X\times X\rightarrow X$ are the two projections.
\end{definition}
The following is a simple but essential result relating the skew Pontrjagin product 
to the surjectivity of the multiplication map (1). It is a restatement of \cite{Pareschi}
(1.1) in a form convenient to us and we reproduce Pareschi's argument for the sake 
of completeness. 

\begin{proposition}\label{41}
Assume that $E\hat{*} F$ is globally generated and that $h^{i}(t_{x}^{*}E
\otimes F)=0$ for all $x\in X$ and all $i>0$, where $t_{x}$ is the translation 
by $x$. Then for all $x\in X$ the multiplication map
$$H^{0}(t_{x}^{*}E)\otimes H^{0}(F)\longrightarrow H^{0}(t_{x}^{*}E\otimes F)$$
is surjective and in particular (1) is surjective. 
\end{proposition}
\begin{proof}
The fact that $h^{i}(t_{x}^{*}E\otimes F)=0$ for all $i>0$ and all $x\in X$ implies 
by base change that $E\hat{*} F$ is locally free with fiber
$$E\hat{*} F(x)\cong H^{0}(t_{x}^{*}E\otimes F).$$
We also have a natural isomorphism
$$\varphi: H^{0}(E\hat{*} F)\overset{\sim}{\longrightarrow}H^{0}(E)\otimes 
H^{0}(F)$$
obtained as follows. One one hand by Leray we naturally have 
$$H^{0}({p_{1}}_{*}((p_{1}+p_{2})^{*}E\otimes p_{2}^{*}F))\cong H^{0}
((p_{1}+p_{2})^{*}E\otimes p_{2}^{*}F)$$
and on the other hand the automorphism $(p_{1}+p_{2}, p_{2})$ of $X\times X$
induces an isomorphism
$$H^{0}((p_{1}+p_{2})^{*}E\otimes p_{2}^{*}F)\cong H^{0}(E)\otimes H^{0}(F),$$
so $\varphi$ is obtained by composition. If we identify $H^{0}(E)\otimes H^{0}(F)$
with both $H^{0}(E\hat{*} F)$ (via $\varphi$) and $H^{0}(t_{x}^{*}E)\otimes
H^{0}(F)$ (via $t_{x}^{*}\times id$), then it is easily seen that the multiplication
map 
$$H^{0}(t_{x}^{*}E)\otimes H^{0}(F)\longrightarrow H^{0}(t_{x}^{*}E\otimes F)$$
coincides with the evaluation map
$$H^{0}(E\hat{*} F)\overset{ev_{x}}{\longrightarrow} E\hat{*} F(x)$$
and this proves the assertion.
\end{proof}

\begin{remark}
There is a clear relationship between the skew Pontrjagin product and the usual 
Pontrjagin product as defined in \ref{30}(2). In fact (see \cite{Pareschi} (1.2)):
$$E\hat{*} F\cong E*(-1_{X})^{*}F.$$
If $F$ is symmetric the two notions coincide and one can hope to apply results 
like \ref{30}(2). This whole circle of ideas will be used in Section 5 below. 
\end{remark}

Finally we extract another result on multiplication maps that will be useful in the 
sequel. It is a particular case of \cite{Pareschi} (3.8), implicit in Kempf's work
\cite{Kempf} on syzygies of abelian varieties. 

\begin{proposition}\label{42}
Let $E$ be a vector bundle on $X$, $L$ an ample line bundle and $m\geq 2$ an integer.
Assume that
$$h^{i}(E\otimes L^{\otimes k}\otimes \alpha)=0,~ \forall i>0,~\forall k\geq -2 {\rm~ and~} 
\forall \alpha\in {\rm Pic}^{0}(X).$$
Then the multiplication map
$$H^{0}(L^{\otimes m})\otimes H^{0}(E)\longrightarrow H^{0}(L^{\otimes m}\otimes E)$$
is surjective.
\end{proposition}

\section{\textbf{The Verlinde bundles} $E_{r,k}$}

In the present section we define the main objects of this paper and study their first 
properties. These are vector bundles on the Jacobian of $X$
which are naturally associated to generalized theta line bundles on the moduli spaces
$U_{X}(r,0)$ and play a key role in what follows.  

\begin{definition} 
For every positive
integers $r$ and $k$ and every $N\in$Pic$^{g-1}(X)$, the $(r,k)$ -$\it{Verlinde}$
bundle on $J(X)$ is the vector bundle $$E_{r,k}(=E^{N}_{r,k})
:={\rm det}_{*}\OO(k\Theta_{N}),$$
where ${\rm det}:  U_{X}(r,0) \rightarrow J(X)$ is the determinant map. The
dependence of the definition on the choice of $N$ is implicitly assumed, but
not emphasized by the notation. Recall that by convention 2.2 we will (and 
it is enough) to assume that $N$ is a theta characteristic. Also, most of the time
the rank $r$ is fixed and we refer to $E_{r,k}$ as the $\it{level~ k}$ 
Verlinde bundle.  
\end{definition}

\begin{remark} The fibers of the determinant map are the moduli spaces of
vector bundles of fixed determinant $SU_{X}(r,L)$ with $L\in$Pic$^{0}(X)$.
The restriction of $\OO(k\Theta_{N})$ to such a fiber is exactly $\LL^{k}$,
where $\LL$ is the determinant bundle. Since on $SU_{X}(r)$ the dualizing 
sheaf is isomorphic to $\LL^{-2r}$, by the rational singularities version of
the Kodaira vanishing theorem 
these have no higher
cohomology and it is clear that $E_{r,k}$ is a vector bundle of rank
$s_{r,k}:=h^{0}(SU_{X}(r),\LL^{k})$. The fiber of $E_{r,k}$ at a point $L$
is naturally isomorphic to the Verlinde vector space
$H^{0}(SU_{X}(r,L),\LL^{k})$ and this justifies the terminology "Verlinde
bundle" that we are using. 
\end{remark} 
 
Much of the study of the vector
bundles $E_{r,k}$ is governed by the fact that they decompose very nicely
when pulled back via multiplication by $r$. To see this, first recall from
\cite{Donagi}, \S2 and \S4, that there is a cartesian diagram:  
$$\xymatrix{
SU_{X}(r)\times J(X) \ar[r]^{\hspace{7mm} \tau} \ar[d]_{p_{2}} & U_{X}(r,0) \ar[d]^{{\rm det}} \\
J(X) \ar[r]^{r_{J}} & J(X) } $$
where $\tau$ is the tensor product of vector bundles, $p_{2}$ is the projection
on the second factor and $r_{J}$ is multiplication by $r$. 
The top and bottom maps are \'etale covers of degree
$r^{2g}$ and one finds in \cite{Donagi} the formula:  
\begin{equation}
\tau^{*}\OO(\Theta_{N})\cong \LL\boxtimes\OO_{J}(r\Theta_{N}).
\end{equation} 
Using the notation $V_{r,k}:=H^{0}(SU_{X}(r),\LL^{k})$, we
have the following simple but very important 

\begin{lemma}\label{1}
$r_{J}^{*}E_{r,k}\cong V_{r,k}\otimes\OO_{J}(kr\Theta_{N}).$ 
\end{lemma}

\begin{proof} By the push-pull formula (see \cite{Hartshorne}, III.9.3) and (2)
we obtain: $$r_{J}^{*}E_{r,k}\cong r_{J}^{*}{\rm det}_{*}\OO(k\Theta_{N})\cong
{p_{2}}_{*}\tau^{*}\OO(k\Theta_{N})$$ $$\cong
{p_{2}}_{*}(\LL^{k}\boxtimes\OO_{J}(kr\Theta_{N}))\cong
V_{r,k}\otimes\OO_{J}(kr\Theta_{N}).$$ 
\end{proof} 
An immediate consequence
of this property is the following (recall that a vector bundle is
called \emph{polystable} if it decomposes as a direct sum of stable
bundles of the same slope): 

\begin{corollary}\label{3} 
$E_{r,k}$ is an ample
vector bundle, polystable with respect to any polarization on $J(X)$. 
\end{corollary} 

\begin{proof} Both properties can be checked up to finite covers (see
e.g. \cite{Lehn} \S3.2) 
and they are obvious for $r_{J}^{*}E_{r,k}$.  
\end{proof} 
For future reference it is also necessary to study how the bundles $E_{r,k}$ behave
under the Fourier-Mukai transform (recall the definitions from 2.3). 

\begin{lemma} 
$E_{r,k}$ satisfies I.T.
with index 0 and so $\widehat{E_{r,k}}$ is a vector bundle of rank
$h^{0}(U_{X}(r,0), \OO(k\Theta_{N}))$ satisfying I.T. with index $g$. 
\end{lemma} 
\begin{proof} 
Using the usual identification between Pic$^{0}(J(X))$ and Pic$^{0}(X)$,
we have to show that $H^{i}(E_{r,k}\otimes P)=0$ for all $P\in$Pic$^{0}(X)$
and all $i>0$. 
But $H^{i}(E_{r,k}\otimes P)$ is a direct summand in
$H^{i}({r_{J}}_{*}r_{J}^{*}(E_{r,k}\otimes P))$ and so it is enough to
have
the vanishing of $H^{i}(r_{J}^{*}(E_{r,k}\otimes P))$. This is obvious by 
the formula in \ref{1}. By \ref{30}(5) we have:
$$rk(\widehat{E_{r,k}})=h^{0}(E_{r,k})=h^{0}(U_{X}(r,0),\OO(k\Theta_{N})).$$
The last statement follows also from \ref{30}(5).  
\end{proof} 
By Serre duality, one obtains in a similar way the statement: 

\begin{lemma}\label{13} 
$E_{r,k}^{*}$ satisfies I.T. with
index $g$ and so $\widehat{E_{r,k}^{*}}$ is a vector bundle of rank
$h^{0}(U_{X}(r,0),\OO(k\Theta_{N}))$ satisfying I.T. with index 0.
\end{lemma}
The Verlinde bundles are particularly easy to compute when $k$ is a multiple of 
the rank $r$, but note that for general $k$ the decomposition of $E_{r,k}$ into stable factors is not 
so easy to describe.

\begin{proposition}\label{12}
For all $m\geq 1$
$$E_{r,mr}\cong 
\underset{s_{r,mr}}{\bigoplus}\OO_{J}(m\Theta_{N})
(\cong V_{r,mr}\otimes\OO_{J}(m\Theta_{N})).$$
\end{proposition}
\begin{proof}
By \ref{1} we have $$r_{J}^{*}E_{r,mr}\cong V_{r,mr}\otimes\OO_{J}(mr^{2}\Theta_{N}).$$ 
Since $N$ is a theta characteristic, $\Theta_{N}$ is symmetric and so 
$r_{J}^{*}\OO(m\Theta_{N})\cong \OO_{J}(mr^{2}\Theta_{N})$. We get 
$$r_{J}^{*}E_{r,mr}\cong r_{J}^{*}(V_{r,mr}\otimes\OO(m\Theta_{N})).$$
Note now that the diagram giving \ref{1} is equivariant with respect to the $X_{r}$, the 
group of $r$-torsion points of $J(X)$, if we let $X_{r}$ act on $SU_{X}(r)\times J(X)$ 
by $(\xi,(F,L))\rightarrow (F,L\otimes \xi)$, by twisting on $U_X(r,0)$, by translation on 
the left $J(X)$ and trivially on the right $J(X)$.
Also the action of $X_{r}$ on $\tau^{*}\OO(mr\Theta_{N})\cong \LL^{mr}\boxtimes 
\OO_{J}(mr^{2}\Theta_{N})$ is on $\OO_{J}(mr^{2}\Theta_{N})$ the same as the natural (pullback)
action. By chasing the diagram we see then that the vector bundle isomorphism above is 
equivariant with respect to the natural $X_{r}$ action on both sides. Since we have an 
an induced isomorphism:
$${r_{J}}_{*}r_{J}^{*}E_{r,mr}\cong {r_{J}}_{*}r_{J}^{*}(V_{r,mr}\otimes 
\OO(m\Theta_{N}))$$ 
and $E_{r,mr}$ and  
$V_{r,mr}\otimes \OO_{J}(m\Theta_{N})$ are both eigenbundles
with respect to the trivial character, the lemma follows.
\end{proof}

The rest of the section will be devoted to a further study of these bundles in 
the case $k=1$,
where more tools are available. The main result is that $E_{r,1}$ is a simple 
vector bundle and
this fact will be exploited in the next section. We show this after proving a 
very simple lemma, to
the effect that twisting by $r$-torsion line bundles does not change $E_{r,1}$.
  
\begin{lemma}\label{2}
$E_{r,1}\otimes P_{\xi}\cong E_{r,1}$ for any $r$-torsion line bundle $P_{\xi}$ 
on $J(X)$ corresponding to an $r$-torsion
$\xi\in {\rm Pic}^{0}(X)$ by the usual identification.  
\end{lemma} 
\begin{proof} 
By definition and the projection formula one
has 
$$E_{r,1}\otimes P_{\xi}\cong {\rm det}_{*}(\OO(\Theta_{N})\otimes {\rm det}^{*}P_{\xi})\cong
{\rm det}_{*}\OO(\Theta_{N\otimes \xi}),$$ 
where the last isomorphism is an application of \ref{32}.
Now the following commutative diagram:  
$$\xymatrix{
U_{X}(r,0) \ar[r]^{\otimes\xi} \ar[d]_{{\rm det}} & U_{X}(r,0) \ar[d]^{{\rm det}} \\
J(X) \ar[r]^{id} & J(X) } $$ 
shows that ${\rm det}_{*}\OO(\Theta_{N})\cong 
{\rm det}_{*}\OO(\Theta_{N\otimes \xi})$, which is exactly the statement of the lemma.
\end{proof}  

\begin{proposition} 
$E_{r,1}$ is a simple vector bundle. 
\end{proposition} 
\begin{proof}
This
follows from a direct computation of the number of endomorphisms
of $E_{r,1}$.  By lemma \ref{1} and
the Verlinde formula at level $1$, $r_{J}^{*}E_{r,1}\cong \underset{r^{g}}{\bigoplus}
\OO_{J}(r\Theta_{N})$. Then: 
$$h^{0}({r_{J}}_{*}r_{J}^{*}(E_{r,1}^{*}\otimes
E_{r,1}))=h^{0}(r_{J}^{*}(E_{r,1}^{*}\otimes
E_{r,1}))$$ $$ =h^{0}((\underset{r^{g}}{\oplus}\OO_{J}(-r\Theta_{N}))\otimes
(\underset{r^{g}}{\oplus}\OO_{J}(r\Theta_{N})))=r^{2g}.$$ 
On the other hand, since $r_{J}$
is a Galois cover with Galois group $X_{r}$,
 we have the formula ${r_{J}}_{*}\OO_{J}\cong 
\underset{\xi\in
X_{r}}{\bigoplus}P_{\xi}$. Combined with lemma \ref{2}, this gives
$${r_{J}}_{*}r_{J}^{*}(E_{r,1}^{*}\otimes E_{r,1})\cong 
\underset{\xi\in X_{r}}{\bigoplus}E_{r,1}^{*}\otimes E_{r,1}\otimes P_{\xi}\cong 
\underset{r^{2g}}{\bigoplus} E_{r,1}^{*}\otimes E_{r,1}.$$
The two relations imply that $h^{0}(E_{r,1}^{*}\otimes
E_{r,1})=1$, so $E_{r,1}$ is simple.
\end{proof}
\begin{remark}
Using an argument similar to the one given above, plus the Verlinde formula, 
it is not hard to see that the bundles $E_{r,k}$ are not simple for $k\geq 2$.
In the case when $k$ is a multiple of $r$ they even decompose as direct 
sums of line bundles, as we have already seen in \ref{12}. This shows why some
special results that we will obtain in the case $k=1$ do not admit 
straightforward extensions to higher $k$'s.
\end{remark}

\section{\textbf{Stability of Fourier transforms and duality for generalized theta functions}}

One of the main features of the vector bundles $E_{r,1}$, already
observed in the previous section, is that they are simple.
We will see below that in fact they are even stable (with respect to any polarization 
on $J(X)$).
This fact, combined with the fact that
$\widehat{\OO_{J}(r\Theta_{N})}$ is also stable (see Proposition \ref{7} below),
gives simple proofs of some results
of Beauville-Narasimhan-Ramanan \cite{BNR} and Donagi-Tu \cite{Donagi}.  
Note that for the proofs of these applications it is enough to use the
simpleness 
of the vector bundles mentioned above. 
 
\begin{proposition}\label{7}
$\widehat{E_{1,r}}=\widehat{\OO_{J}(r\Theta_{N})}$ is stable with
respect to any polarization on $J(X)$.
\end{proposition}

This is a consequence of a more
general fact of independent interest, saying that this is indeed true for an arbitrary
nondegenerate line bundle on an abelian variety. Recall that a line bundle $A$ on the 
abelian variety $X$ is called $\it{nondegenerate}$ if $\chi(A)\neq 0$. By \cite{Mumford2} \S16
this implies that there is a unique $i$ (the $\it{index}$ of $A$) such that $H^{i}(A)\neq 0$.  

\begin{proposition}
Let $A$ be a nondegenerate line bundle on an abelian variety $X$. The
Fourier-Mukai transform $\widehat{A}$ is stable with respect
to any polarization on $X$. 
\end{proposition}
\begin{proof} 
Let's begin by fixing a polarization on $\widehat{X}$, so that
stability will be understood with respect to this polarization. 
Consider the isogeny defined by $A$
$$\begin{array}{cccc}
\phi_{A}:& X& \longrightarrow & {\rm Pic}^{0}(X)\cong\widehat{X}\\
&x &\leadsto& t_{x}^*A\otimes A^{-1}
\end{array}
$$
If $i$ is the index of $A$, it follows from \ref{30}(3) that
$\phi_{A}^{*}\widehat{A}\cong V\otimes A^{-1}$, where
$V:=H^{i}(A)$. As we already mentioned in the proof of \ref{3}, by \cite{Lehn} \S3.2 
this already implies that $\widehat{A}$ is polystable. On the other hand, by 
\ref{30}(5) the Fourier 
transform of any line bundle is simple, so $\widehat{A}$ must be stable.
\end{proof} 

\begin{remark}
It is worth noting that one can avoid the use of \cite{Lehn} \S3.2 quoted above 
and use only the easier fact that semistability is preserved by finite covers. 
More precisely, $\widehat{A}$ has to be semistable, but assume that it is not stable. 
Then we can choose a
maximal destabilizing subbundle $F\subsetneq \widehat{A}$, which
must obviously be semistable and satisfy
$\mu(F)=\mu(\widehat{A})$. Again $\phi_{A}^{*}F$ must be semistable, with respect
 to the 
pull-back polarization. But $\phi_{A}^{*}F\subsetneq
\phi_{A}^{*}\widehat{A}\cong V\otimes A^{-1}$ and by semistability
this implies that $\phi_{A}^{*}F\cong V^{\prime}\otimes A^{-1}$ with
$V^{\prime}\subsetneq V$.
This situation is overruled by the presence of the action of Mumford's
theta-group $\mathcal{G}(A)$. Recall from \cite{Mumford} that $A$ is endowed with 
a natural
$\mathcal{G}(A)$-linearization of weight $1$. On the other hand,
$\phi_{A}^{*}F$ has a natural $K(A)$-linearization which can be
seen as a $\mathcal{G}(A)$-linearization of weight 0, where $K(A)$ is the kernel of 
the isogeny $\phi_{A}$ above. By
tensoring we obtain a weight $1$ $\mathcal{G}(A)$-linearization on 
$V^{\prime}\otimes\OO_{X}\cong \phi_{A}^{*}F\otimes A$ and thus an
induced weight 1 representation on $V^{\prime}$. It is known
though from \cite{Mumford} that $V$ is the unique irreducible
representation of $\mathcal{G}(A)$ up to isomorphism, so we get a contradiction. 
\end{remark}

We now return to the study of the relationship between the
bundles $E_{r,k}$ and their Fourier transforms. 
First note that throughout the rest of the paper we will use the fact that 
$J(X)$ is canonically isomorphic to its dual ${\rm Pic}^{0}(J(X))$ and consequently 
we will use the same notation for both. Thus all the Fourier transforms should 
be thought of as coming from the dual via this identification, although since
there is no danger of confusion this will not be visible in the notation. 
The fiber of $E_{r,k}$ over a point $\xi\in J(X)$ is $H^{0}(SU_{X}(r,\xi),\LL^{k})$, 
while the fiber of $\widehat{E_{k,r}}$ over the same point is canonically 
isomorphic to $H^{0}(J(X),E_{k,r}\otimes P_{\xi})$. By \ref{32} the latter is 
isomorphic to $H^{0}(U_{X}(k,0),\OO(r\Theta_{N\otimes \eta}))$, where 
$\eta^{\otimes r}\cong \xi$ (it is easy to see that this does not depend on the 
choice of $\eta$). The strange duality conjecture (see \cite{Beauville} \S8 and 
\cite{Donagi} \S5) says that there is a canonical isomorphism:
$$H^{0}(SU_{X}(r,\xi),\LL^{k})^{*}\cong H^{0}(U_{X}(k,0),
\OO(r\Theta_{N\otimes \eta}))$$ which will be described more precisely later 
in this section. This suggests then that one should relate somehow the vector 
bundles $E_{r,k}$ and $\widehat{E_{k,r}}$ via the diagram:
$$\xymatrix{
U_{X}(r,0) \ar[dd]_{{\rm det}} & & U_{X}(k,0) \ar[dd]^{{\rm det}} \\
&  J(X)\times J(X) \ar[dl]_{p_{1}} \ar[dr]^{p_{2}} \\
J(X) & & J(X) }$$
The first proposition treats the case $k=1$ and 
establishes the fact that the dual of $E_{r,1}$ is
nothing else but the Fourier transform of $\OO_{J}(r\Theta_{N})$. 

\begin{proposition}\label{6}
$E_{r,1}^{*}\cong \widehat{E_{1,r}}.$
\end{proposition}
\begin{proof}
By the duality theorem \ref{34}, it is enough to
show that $\widehat{E_{r,1}^{*}}\cong (-1_{J})^{*}E_{1,r}$. But
$E_{1,r}$ is just $\OO_{J}(r\Theta_{N})$, which is symmetric
since $N$ is a theta characteristic. So what we have to prove is 
$$\widehat{E_{r,1}^{*}}\cong \OO_{J}(r\Theta_{N}).$$
Since $r_{J}$ is Galois with Galois group $X_{r}$, we have (see e.g \cite{Narasimhan}
(2.1)): 
$$r_{J}^{*}{r_{J}}_{*}\widehat{E_{r,1}^{*}}\cong \underset{\xi\in X_{r}}
{\bigoplus}t_{\xi}^{*}\widehat{E_{r,1}^{*}}.$$
But translates commute with tensor products via the Fourier transform (cf. \ref{30}(4))
and so 
$$t_{\xi}^{*}\widehat{E_{r,1}^{*}}\cong \widehat{E_{r,1}^{*}\otimes P_{\xi}}
\cong E_{r,1}^{*},$$
where as usual $P_{\xi}$ is the line bundle in Pic$^{0}(J(X))$ that corresponds 
to $\xi$ and the last isomorphism follows from \ref{2}. Thus we get 
$$r_{J}^{*}{r_{J}}_{*}\widehat{E_{r,1}^{*}}\cong \underset{r^{2g}}
{\bigoplus}\widehat{E_{r,1}^{*}}$$
and the idea is to compute this bundle in a different way, by using the 
behavior of the Fourier transform under isogenies. More precisely, by applying 
\ref{30}(1) we get the isomorphisms:
$$r_{J}^{*}{r_{J}}_{*}\widehat{E_{r,1}^{*}}\cong r_{J}^{*}
\widehat{r_{J}^{*}E_{r,1}^{*}}\cong r_{J}^{*}(V_{r,1}^{*}\otimes 
\OO_{J}(-r\Theta_{N}))^{\widehat{}}~$$ 
$$\cong \underset{r^{g}}{\bigoplus}r_{J}^{*}\widehat{\OO_{J}(-r\Theta_{N})}
\cong \underset{r^{2g}}{\bigoplus}\OO_{J}(r\Theta_{N}).$$
The second isomorphism follows from \ref{1}, while the fourth follows from 
\ref{30}(3) and the Verlinde formula at level $1$. The outcome is the 
isomorphism $$\underset{r^{2g}}{\bigoplus}\widehat{E_{r,1}^{*}}\cong \underset{r^{2g}}{\bigoplus}\OO_{J}(r\Theta_{N}).$$ Now $\widehat{E_{r,1}^{*}}$ is simple 
since $E_{r,1}$ is simple, so the previous isomorphism implies the stronger
fact that $$\widehat{E_{r,1}^{*}}\cong \OO_{J}(r\Theta_{N}).$$  
\end{proof}

\begin{example}
The relationship between the Chern character of a sheaf and that of its 
Fourier transform established in \cite{Mukai2} (1.18) allows us to easily 
compute the first Chern class of $E_{r,1}$ as a consequence of the previous
proposition. More precisely, if 
$\theta$ is the class of a theta divisor on $J(X)$, we have 
$$c_{1}(\widehat{\OO_{J}(r\Theta_{N})})={\rm ch}_{1}(\widehat{\OO_{J}(r\Theta_{N})})
=(-1)\cdot PD_{2g-2}({\rm ch}_{g-1}(\OO_{J}(r\Theta_{N})))=$$ 
$$=(-1)\cdot PD_{2g-2}(r^{g-1}\cdot\theta^{g-1}/(g-1)!)= -r^{g-1}\cdot\theta.$$
where $PD_{2g-2}:H^{2g-2}(J(X),\ZZ)\rightarrow H^{2}(J(X),\ZZ)$ is the 
Poincar\'e duality. We get: 
$$c_{1}(E_{r,1})=r^{g-1}\cdot\theta.$$   
\end{example}

We are able to prove the analogous fact for higher $k$'s only modulo multiplication
by $r$:

\begin{proposition}\label{35}
$r_{J}^{*}E_{r,k}^{*}\cong r_{J}^{*}\widehat{E_{k,r}}.$
\end{proposition}
\begin{proof}
We know that $k_{J}^{*}E_{k,r}\cong
\underset{s_{k,r}}{\bigoplus}\OO_{J}(kr\Theta_{N})$, so by \ref{30}(1) 
we obtain 
$${k_{J}}_{*}\widehat{E_{k,r}}\cong
\underset{s_{k,r}}{\bigoplus}\widehat{\OO_{J}(kr\Theta_{N})}.$$
As in the previous proposition, since the Galois group of $k_{J}$ is 
$X_{k}$, we have
$$k_{J}^{*}{k_{J}}_{*}\widehat{E_{k,r}}\cong \underset{\xi\in X_{k}}
{\bigoplus}t_{\xi}^{*}\widehat{E_{k,r}} ~{\rm
and~so}~
(kr)_{J}^{*}{k_{J}}_{*}\widehat{E_{k,r}}\cong \underset{\xi\in X_{k}}
{\bigoplus}r_{J}^{*}t_{\xi}^{*}\widehat{E_{k,r}}.$$
Moreover the isomorphisms above show as before that $\widehat{E_{k,r}}$ has to be 
semistable with respect to an arbitrary polarization. On the other 
by \ref{30}(3) we have
$$(kr)_{J}^{*}{k_{J}}_{*}\widehat{E_{r,k}}\cong
\underset{s_{k,r}}{\bigoplus}(kr)_{J}^{*}\widehat{\OO_{J}(kr\Theta_{N})}
\cong
\underset{k^{g}r^{g}s_{k,r}}{\bigoplus}\OO_{J}(-kr\Theta_{N}).$$ 
This gives us the isomorphism $$\underset{\xi\in X_{k}}{\bigoplus}r_{J}^{*}
t_{\xi}^{*}\widehat{E_{k,r}}
\cong\underset{k^{g}r^{g}s_{k,r}}{\bigoplus}\OO_{J}(-kr\Theta_{N})$$ which
in particular implies (recall semistability) that  
$$r_{J}^{*}\widehat{E_{k,r}}\cong \underset{s_{r,k}}{\bigoplus}\OO_{J}(-kr\Theta_{N})
\cong r_{J}^{*}E_{r,k}^{*}.$$
The important fact that the last index of summation is $s_{r,k}$ follows from the well-known
``symmetry'' of the Verlinde formula, which is:
$$r^{g}s_{k,r}=k^{g}s_{r,k}.$$
\end{proof}

The propositions above allow us to give alternative proofs of some results in
\cite{BNR} and \cite{Donagi} concerning duality between spaces of 
generalized theta functions. We first show how one can recapture a theorem 
of Donagi-Tu in the present context. The full version of the theorem (i.e. for 
arbitrary degree) can be obtained by the same method (cf. Section 6).

\begin{theorem}\label{4}(\cite{Donagi}, Theorem 1)
For any $L\in {\rm Pic}^{0}(X)$ and any $N\in {\rm Pic}^{g-1}(X)$, we have:
$$h^{0}(SU_{X}(k,L),\LL^{r})\cdot r^{g}=h^{0}(U_{X}(k,0),\OO(r\Theta_{N}))
\cdot k^{g}.$$
\end{theorem}
\begin{proof}
We will actually prove the following equality:
$$h^{0}(SU_{X}(r,L),\LL^{k})=h^{0}(U_{X}(k,0),\OO(r\Theta_{N})).$$
The statement will then follow from the same symmetry of the Verlinde 
formula $r^{g}s_{k,r}=k^{g}s_{r,k}$ mentioned in the proof of \ref{35}. 
To this end we can use Proposition \ref{35} to obtain ${\rm rk}(E_{r,k}^{*})=
{\rm rk}(\widehat{E_{k,r}})$. But on one hand 
$${\rm rk}(E_{r,k}^{*})=h^{0}(SU_{X}(r,L),\LL^{k})$$
while on the other hand by \ref{30}(5)
$${\rm rk}(\widehat{E_{k,r}})=h^{0}(J(X),E_{k,r})=h^{0}(U_{X}(k,0),\OO(r\Theta_{N}))$$
as required.
\end{proof}

It is worth mentioning that since we are assuming the Verlinde formula all 
throughout, an important particular case of the theorem above is:

\begin{corollary}(\cite{BNR}, Theorem 2)
$h^{0}(U_{X}(r,0),\OO(\Theta_{N}))=1.$
\end{corollary}

As it is very well known, this fact is essential in setting up the strange 
duality. Turning to an application in this direction, the results above also 
lead to a simple proof of the strange duality at level $1$, which has been 
first given a proof in \cite{BNR}, Theorem 3. It is important though to 
emphasize again that here, unlike in the quoted paper, the Verlinde formula is granted.
So the purpose of the next application is to show how, with the knowledge
of the Verlinde numbers, the strange duality at level $1$ can simply be seen as the 
solution of a stability problem for vector bundles on $J(X)$.
This may also provide a global method for understanding the conjecture for higher 
levels. A few facts in this direction will be mentioned at the end of this 
section (cf. Remark \ref{51}).

Before turning to the proof, we need the following general result, 
which is a globalization of \S3 in \cite{Donagi}.

\begin{proposition}\label{36}
Consider the tensor product map $$\tau: U_{X}(r,0)\times
U_{X}(k,0) \longrightarrow U_{X}(kr,0)$$ and the map
$$\phi:={\rm det}\times {\rm det}: U_{X}(r,0)\times
U_{X}(k,0) \longrightarrow J(X)\times J(X).$$ Then
$$\tau^{*}\OO(\Theta_{N})\cong
p_{1}^{*}\OO_{J}(k\Theta_{N})\otimes
p_{2}^{*}\OO_{J}(r\Theta_{N})\otimes \phi^{*}\mathcal{P},$$ where
$\mathcal{P}$ is a Poincar\'e line bundle on $J(X)\times J(X)$,
normalized such that $\mathcal{P}_{|\{0\}\times J(X)}\cong \OO_{J}$. 
\end{proposition}
\begin{proof}
We will compare the restrictions of the two line bundles to 
fibers of the projections. First fix $F\in U_{X}(k,0)$. We
have $$\tau^{*}\OO(\Theta_{N})|_{U_{X}(r,0)\times \{F\}}\cong
\tau_{F}^{*}\OO(\Theta_{N})\cong \OO(\Theta_{F\otimes N}),$$
where $\tau_{F}$ is the map given by twisting with $F$. But by \ref{32}
one has 
$$\OO(\Theta_{F\otimes N})\cong \OO(k\Theta_{N})\otimes
{\rm det}^{*}({\rm detF})(\cong \OO(k\Theta_{N})\otimes
{\rm det}^{*}(\mathcal{P}_{|J(X)\times \{{\rm det}F\}})).$$ 
On the other hand obviously 
$${p_{1}^{*}\OO_{J}(k\Theta_{N})\otimes
p_{2}^{*}\OO_{J}(r\Theta_{N})\otimes
\phi^{*}\mathcal{P}}|_{U_{X}(r,0)\times \{F\}}\cong \OO(k\Theta_{N})\otimes
{\rm det}^{*}({\rm det}F).$$ 

Let's now fix $E\in U_{X}(r,0)$ such that ${\rm det}E\cong
\OO_{X}$. Using the same \ref{32} we get
$$\tau^{*}\OO(\Theta_{N})|_{\{E\}\times 
U_{X}(k,0)}\cong
\OO(\Theta_{E\otimes N})\cong \OO(r\Theta_{N})$$ and we also have 
$${p_{1}^{*}\OO_{J}(k\Theta_{N})\otimes
p_{2}^{*}\OO_{J}(r\Theta_{N})\otimes
\phi^{*}\mathcal{P}}|_{\{E\}\times U_{X}(k,0)}$$
$$\cong
\OO(r\Theta_{N})\otimes
{\rm det}^{*}(\mathcal{P}_{|\{\OO_{X}\}\times J(X)})\cong \OO(r\Theta_{N}).$$ 
The desired isomorphism follows now from the see-saw
principle(see e.g. \cite{Mumford2}, I.5.6).
\end{proof}   
Theorem \ref{4} tells us that there is essentially a unique
nonzero section $$s\in H^{0}(U_{X}(kr,0),\OO(\Theta_{N})),$$ which
induces via $\tau$ a nonzero section $$t\in H^{0}(U_{X}(r,0)\times
U_{X}(k,0),\tau^{*}\OO(\Theta_{N})$$ $$\cong H^{0}(U_{X}(r,0)\times
U_{X}(k,0), p_{1}^{*}\OO_{J}(k\Theta_{N})\otimes
p_{2}^{*}\OO_{J}(r\Theta_{N})\otimes \phi^{*}\mathcal{P}).$$ But
notice that from the projection formula we get 
$$\phi_{*}(p_{1}^{*}\OO_{J}(k\Theta_{N})\otimes
p_{2}^{*}\OO_{J}(r\Theta_{N})\otimes \phi^{*}\mathcal{P})$$ 
$$\cong \phi_{*}(p_{1}^{*}\OO_{J}(k\Theta_{N})\otimes
p_{2}^{*}\OO_{J}(r\Theta_{N}))\otimes \mathcal{P}\cong
p_{1}^{*}E_{r,k}\otimes p_{2}^{*}E_{k,r}\otimes \mathcal{P},$$ 
so $t$ induces a section (denoted also by $t$): 
$$0\neq t\in H^{0}(J(X)\times J(X), p_{1}^{*}E_{r,k}\otimes
p_{2}^{*}E_{k,r}\otimes \mathcal{P})\cong H^{0}(E_{r,k}\otimes \widehat{E_{k,r}}).$$
This is nothing else but a globalization of the section
defining the strange duality morphism, as explained in
\cite{Donagi} \S5 (simply because for any $\xi\in$Pic$^{0}(X)$ 
the restriction of $\tau$ in \ref{36} to $SU_{X}(r,\xi)\times U_{X}(k,0)$ is
again the tensor product map 
$$\tau:SU_{X}(r,\xi)\times
U_{X}(k,0) \longrightarrow 
U_{X}(kr,0)$$ and
$\tau^{*}\OO(\Theta_{N})\cong \LL^{k}\boxtimes \OO(r\Theta_{N})$).
In other words, $t$ corresponds to a morphism of vector bundles
$$SD:E_{r,k}^{*}\longrightarrow \widehat{E_{k,r}}$$ which
fiberwise is exactly the strange duality morphism
$$H^{0}(SU_{X}(r,\xi),\LL^{k})^{*}\overset{SD}{\longrightarrow}
H^{0}(U_{X}(k,0),\OO(r\Theta_{N\otimes \eta})),$$
where $\eta^{\otimes r}\cong \xi$.
Since this global morphism 
collects together all the strange duality morphisms as we vary
$\xi$, the strange duality conjecture is equivalent to $SD$
being an isomorphism. 

\begin{conjecture}\label{52}
$SD:E_{r,k}^{*}\longrightarrow \widehat{E_{k,r}}$ is an
isomorphism of vector bundles.
\end{conjecture}
In this context Proposition \ref{35} can be seen as a weak form of ``global'' evidence 
for the conjecture in the case $k\geq 2$. For $k=1$ it can now be easily proved.

\begin{theorem}\label{5}
$SD:E_{r,1}^{*}\longrightarrow \widehat{E_{1,r}}$ is an
isomorphism.
\end{theorem}

\begin{corollary}(cf. \cite{BNR}, Theorem 3)\label{11}
The level 1 strange duality morphism 
$$H^{0}(SU_{X}(r),\LL)^{*}\overset{SD}{\longrightarrow}
H^{0}(J(X),\OO(r\Theta_{N}))$$ is an isomorphism.
\end{corollary}
\begin{proof}(of \ref{5})
All the ingredients necessary for proving this have been
discussed above: by \ref{7} and \ref{6}, the bundles $E_{r,1}^{*}$ and 
$\widehat{E_{1,r}}$
are isomorphic and stable. This means
that $SD$ is essentially the unique nonzero morphism between them,
and it must be an isomorphism.  
\end{proof}

\begin{remark}\label{51}
One step towards \ref{52} is a better understanding of the properties 
of the kernel $F$ of $SD$. A couple of interesting remarks in this direction 
can already be made. Since $E_{r,k}$ and $\widehat{E_{k,r}}$ are polystable of 
the same slope, the same will be true about $F$. On the other hand some simple 
calculus involving \ref{32} and \ref{30}(4) shows that $F$ gets multiplied by a line bundle
in ${\rm Pic}^{0}(J(X))$ when we translate it, so in the language of Mukai (e.g. 
\cite{Mukai2} \S3) it is a semi-homogeneous vector bundle (although clearly not 
homogeneous, i.e. not fixed by all translations).
\end{remark}

\section{\textbf{A generalization of Raynaud's examples}}

In this section we would like to discuss a generalization of the examples
of base points of the theta linear system $|\LL|$ on $SU_{X}(m)$ constructed  
by Raynaud in \cite{Raynaud}. For a survey of this circle of ideas the
reader can consult \cite{Beauville}, \S2. Let us recall here (see \cite{Popa} \S2)
only that for a semistable vector bundle $E$ of rank $m$ to induce a base 
point of $|\LL|$ it is sufficient that it satisfies the property
$$0\leq\mu(E)\leq g-1 {\rm ~and~} h^{0}(E\otimes L)\neq 0 {\rm ~for~}
L\in {\rm Pic}^{0}(X) {\rm ~generic}.$$
The examples of Raynaud are essentially the restrictions of
$\widehat{E_{1,r}^{*}}=\widehat{\OO_{J}(-r\Theta_{N})}$ to some embedding
of the curve $X$ in the Jacobian. We will generalize this by considering 
the Fourier transform of higher level Verlinde bundles 
$\widehat{E_{k,r}^{*}}$ with $k\geq 2$.  

For simplicity, let's fix $k$ and $r$ and denote $F:= 
\widehat{E_{k,r}^{*}}$. This is a vector bundle by lemma \ref{13}. 
Consider also an arbitrary embedding $$j: X
\rightarrow J(X)$$ and denote by $E$ the restriction $F_{|X}$.

\begin{proposition}\label{37} 
E is a semistable vector bundle.
\end{proposition}
\begin{proof}
This follows basically from the proof of Proposition \ref{35}. One can see
in a completely analogous way that:
$$r_{J}^{*}F\cong r_{J}^{*}E_{r,k}\cong \underset{s_{r,k}}
{\bigoplus}\OO_{J}(kr\Theta_{N}).$$
Now if we consider $Y$ to be the preimage of $X$ by $r_{J}$, this shows 
that $r_{J}^{*}F_{|Y}$ is semistable and so $F_{|X}$ is semistable. 
\end{proof}

\begin{proposition} There is an embedding of $X$ in $J(X)$ such that
$E=F_{|X}$ satisfies
$H^{0}(E\otimes L) \neq 0$ for $L\in {\rm Pic}^{0}(X)$ generic.
\end{proposition}

\begin{proof}
The proof goes like in \cite{Raynaud} (3.1) and we repeat it here for
convenience: choose $U$ a nonempty open
subset of $J(X)$ on which $(-1_{J})^{*}E_{k,r}$ is trivial. By Mukai's 
duality theorem \cite{Mukai} (2.2) we know that $\widehat{F} \cong
(-1_{J})^{*}E_{k,r}^{*}$, so there exists a nonzero section $s\in
\Gamma({p_{2}}^{-1}(U), p_{1}^{*}F \otimes \mathcal{P})$. Choose now an
$x_{0}\in J(X)$ such that $s_{|\{x_{0}\}\times U} \neq 0$ and consider an
embedding 
of $X$ in $J(X)$ passing through $x_{0}$. The image of $s$ in
$\Gamma(X\times U,
p_{1}^{*}F \otimes \mathcal{P})$ is nonzero and this implies that
$H^{0}(E \otimes L) \neq 0$ for $L$ generic. 
\end{proof}
We are only left with computing the invariants of $E$.

\begin{proposition}
The rank and the slope of $E$ are given by:
$${\rm rk}(E)=s_{r,k}=h^{0}(SU_{X}(r),\mathcal{L}^{k}) ~{\rm and}~
\mu(E)=gk/r.$$
\end{proposition}
\begin{proof}
By \ref{30}(5) we have
$${\rm rk}(F)={\rm rk}(\widehat{E_{k,r}^{*}})=h^{0}(E_{k,r})=h^{0}(U_{X}(k,0),
\OO(r\Theta_{N})).$$
But from the proof of \ref{4} we know that $h^{0}(U_{X}(k,0),
\OO(r\Theta_{N}))=h^{0}(SU_{X}(r,0),\mathcal{L}^{k})$ and so
${\rm rk}(F)=s_{r,k}$. 
To compute the slope of $E$, first notice that by the proof of Proposition 
\ref{37} we know that $$r^{2g}\cdot\mu(E)={\rm deg}(\OO_{J}(kr\Theta_{N})_{|Y}).$$
But $r_{J}^{*}\OO_{J}(kr\Theta_{N}) \equiv \OO_{J}(kr^{3}\Theta_{N})$, so 
$$r^{2}\cdot
{\rm deg}(\OO_{J}(kr\Theta_{N})_{|Y})={\rm deg}(r_{J}^{*}\OO_{J}(kr\Theta_{N})_{|Y})=
r^{2g}\cdot {\rm deg}(\OO_{J}(kr\Theta_{N})_{|X})=r^{2g+1}kg.$$
Combining the two equalities we get $$\mu(E)=
{\rm deg}(\OO_{J}(kr\Theta_{N})_{|Y})/r^{2g}=gk/r.$$

\end{proof}
In conclusion, for each $r$ and $k$ we obtain a semistable vector
bundle $E$ on $X$ of rank equal to the Verlinde number $s_{r,k}$
and of slope $gk/r$, satisfying the property that $H^{0}(E\otimes
L) \neq 0$ for $L$ generic in Pic$^{0}(X)$. So as long as $k<r$
and $r$ divides $gk$ we obtain new examples of base points for
$|\mathcal{L}|$ on the moduli spaces $SU_{X}(s_{r,k})$, as explained above. 
Raynaud's examples correspond to the case $k=1$. 
See also \cite{Popa} for a study of a
different kind of examples of such base points and for bounds on
the dimension of the base locus of $|\mathcal{L}|$.

\section{\textbf{Linear series on} $U_{X}(r,0)$}

The main application of the Verlinde vector bundles concerns the global generation and 
normal generation of line bundles on $U_{X}(r,0)$. The specific goal is to give effective bounds 
for multiples of the generalized theta line bundles that satisfy the properties mentioned above.
In this direction analogous results for $SU_{X}(r)$ will be used. The starting point is the following 
general result: 
\begin{proposition}\label{8}
Let $f:X\rightarrow Y$ be a flat morphism of projective schemes, with reduced 
fibers. Let $L$ be a line bundle on $X$ and denote $E:= f_{*} L$. Assume that if $X_{y}$ denotes the 
fiber of $f$ over $y\in Y$ the following conditions hold:
\newline
(i) $h^{1}(L)=0$
\newline
(ii) $h^{i}(L_{|X_{y}})=0, \forall y\in Y, \forall i>0$.
\newline
Then $L$ is globally generated as long as $L_{|X_{y}}$ is globally generated 
for all $y\in Y$ and $E$ is globally generated.
\end{proposition}
\begin{proof}
Start with $x\in X$ and consider $y=f(x)$. By $(i)$ we have the exact 
sequence on $X$: 
$$0\longrightarrow H^{0}(L\otimes \II_{X_{y}})\longrightarrow H^{0}(L)
\longrightarrow H^{0}(L_{|X_{y}})\longrightarrow H^{1}(L\otimes \II_{X_{y}})
\longrightarrow 0.$$
The global generation of $L_{|X_{y}}$ implies that there exists a section 
$s\in H^{0}(L_{|X_{y}})$ such that $s(x)\neq 0$. We would like to lift $s$ to 
some  $\bar{s}$, so it is enough to prove that $H^{1}(L\otimes \II_{X_{y}})=0$. The fibers of $f$ are reduced, so $\II_{X_{y}}\cong f^{*}\II_{\{y\}}$. Condition (ii) implies, by the base change theorem, that $R^{i}f_{*}L=0$ for all $i>0$ and 
so by the projection formula we also get $R^{i}f_{*}(L\otimes \II_{X_{y}})=0$ for all $i>0$. The Leray spectral sequence then gives $H^{i}(E)\cong H^{i}(L)$ and 
$H^{i}(E\otimes \II_{\{y\}})\cong H^{i}(L\otimes \II_{X{y}})$ for all $i>0$.
The first isomorphism implies that there is an exact sequence:
$$0\longrightarrow H^{0}(E\otimes \II_{\{y\}})\longrightarrow H^{0}(E) \overset{ev_{y}}{\longrightarrow} H^{0}(E_{y})\longrightarrow H^{1}(E\otimes \II_{\{y\}})\longrightarrow 0.$$
But $E$ is globally generated, which means that $ev_{y}$ is surjective. This implies the vanishing of $H^{1}(E\otimes \II_{\{y\}})$, which by the second isomorphism is 
equivalent to $H^{1}(L\otimes \II_{X_{y}})=0$.  
\end{proof}

The idea is to apply this result to the situation when the map is 
${\rm det}:U_{X}(r,0)\rightarrow J(X)$, $L=\OO(k\Theta_{N})$ and $E=E_{r,k}$. 
Modulo detecting for what $k$
global generation is attained, the conditions of the proposition are satisfied.
The fiberwise global generation problem (i.e. the $SU_{X}(r)$  case) has been given some 
effective solutions in the literature. 
The most recent is the author's result in \cite{Popa2}, improving earlier bounds of Le Potier
\cite{Le Potier} and Hein \cite{Hein}.
We show in \cite{Popa2} \S4 that $\LL^{k}$ on $SU_{X}(r)$ is globally generated 
if $k\geq \frac{(r+1)^{2}}{4}$. We now turn to the effective statement for $E_{r,k}$ and in this direction 
we make essential use of Pareschi's cohomological criterion described in Section 1. 

\begin{proposition}\label{20}
$E_{r,k}$ is globally generated if and only if $k\geq r+1$. 
\end{proposition}
\begin{proof}
The trick is to write $E_{r,k}$ as $E_{r,k}\otimes \OO_{J}(-\Theta_{N})\otimes \OO_{J}(\Theta_{N})$. 
Denote $E_{r,k}\otimes \OO_{J}(-\Theta_{N})$ by $F$. Pareschi's criterion \ref{40} says 
in our case that $E_{r,k}$ is globally generated as long as the condition
$$h^{i}(F\otimes \alpha)=0,~ \forall\alpha\in{\rm Pic}^{0}(X),~\forall i>0$$ 
is satisfied (in fact under this assumption $F\otimes A$ will be globally generated for every 
ample line bundle $A$). 
Arguing as usual, $h^{i}(F\otimes\alpha)=0$ is implied by 
$h^{i}(r_{J}^{*}(F\otimes\alpha))=0$. We have 
$$r_{J}^{*}F\cong r_{J}^{*}E_{r,k}\otimes r_{J}^{*}\OO_{J}(-\Theta_{N})
\cong \OO_{J}((kr-r^{2})\Theta_{N})$$
and this easily gives the desired vanishing for $k\geq r+1$. On the other hand from \ref{12} we know that 
$E_{r,r}\cong \underset{s_{r,r}}{\bigoplus}\OO_{J}(\Theta_{N})$, which is clearly not globally generated.
This shows that the bound is optimal.
\end{proof}
Combining all these we obtain effective bounds for global generation on $U_{X}(r,0)$ in terms of the 
analogous bounds on $SU_{X}(r)$. We prefer to state the general result in a non-effective form though, in order 
to emphasize that it applies algorithmically (but see the corollaries for effective statements):

\begin{theorem}\label{44}
$\OO(k\Theta_{N})$ is globally generated on $U_{X}(r,0)$ as long as $k\geq r+1$
and $\LL^{k}$ is globally generated on $SU_{X}(r)$. Moreover, $\OO(r\Theta_{N})$ is 
not globally generated.
\end{theorem}
\begin{proof}
The first part follows by puting together \ref{20} and \ref{8} in our particular setting. 
To prove that $\OO(r\Theta_{N})$ is not globally generated, let us 
begin by assuming the contrary.
Then the restriction $\LL^{r}$ of $\OO(r\Theta_{N})$ to any of the fibers
$SU_{X}(r,L)$ is also globally generated.

Choose in particular a line bundle $L$ in the support of $\Theta_{N}$ 
on $J(X)$. Restriction to the fiber gives the following long exact sequence
on cohomology: 
$$0\longrightarrow H^{0}(\OO(r\Theta_{N})\otimes \II_{SU_{X}(r,L)})
\longrightarrow H^{0}(\OO(r\Theta_{N}))\overset{\alpha}{\longrightarrow}$$ 
$$\longrightarrow H^{0}(\LL^{r})\longrightarrow H^{1}(\OO(r\Theta_{N})\otimes 
\II_{SU_{X}(r,L)})
\longrightarrow 0.$$ As in the proof of \ref{8}, this sequence can be written
in terms of the cohomology of $E_{r,r}$:
$$0\longrightarrow H^{0}(E_{r,r}\otimes \II_{\{L\}})
\longrightarrow H^{0}(\OO(r\Theta_{N})) 
\overset{\alpha}{\longrightarrow} H^{0}(\LL^{r})\longrightarrow 
H^{1}(E_{r,r}\otimes \II_{\{L\}})\longrightarrow 0.$$
The assumption on $\OO(r\Theta_{N})$ ensures the fact that the map $\alpha$ 
in the sequence above is nonzero, and as a result $$h^{1}(E_{r,r}\otimes 
\II_{\{L\}})< h^{0}(\LL^{r}) = s_{r,r}.$$ On the other hand we use again the fact \ref{12} 
that $E_{r,r}$ is isomorphic to
$\underset {s_{r,r}}{\bigoplus}\OO_{J}(\Theta_{N})$. The additional hypothesis that 
$L\in \Theta_{N}$ says then that $$h^{1}(E_{r,r}\otimes \II_{\{L\}})= 
s_{r,r}$$ which is a contradiction.
\end{proof}
As suggested above, by combining this with the effective result on $SU_{X}(r)$ given in 
\cite{Popa2} \S4, we get: 

\begin{corollary}(cf. \cite{Popa2} (5.3))\label{22}
$\OO(k\Theta_{N})$ is globally generated on $U_{X}(r,0)$ for  
$$k\geq \frac{(r+1)^{2}}{4}.$$
\end{corollary}
It is important, as noted in the introduction, to emphasize the fact that the $SU_{X}(r)$ bound 
may still allow for improvement. Thus the content and formulation of the theorem 
certainly go beyond this corollary.
For moduli spaces of vector bundles of rank $2$ and $3$ though, in view of the second part of \ref{44}, 
we actually have optimal results (see also \cite{Popa2} (5.4)).

\begin{corollary}\label{60}
(i) $\OO(3\Theta_{N})$ is globally generated on $U_{X}(2,0)$.
\newline
(ii) $\OO(4\Theta_{N})$ is globally generated on $U_{X}(3,0)$. 
\end{corollary}
These are natural extensions of the fact that $\OO(2\Theta_{N})$
is globally generated on $J(X)\cong U_{X}(1,0)$ (see e.g. \cite{Griffiths}, Theorem 2, p.317). In 
\cite{Popa2} \S5 we also state some questions and conjectures about optimal bounds in general.

\begin{remark}
A similar technique can be applied to study the base point freeness 
of more general linear series on $U_{X}(r,0)$. 
This is done in \cite{Popa2} (5.9).
\end{remark}

\begin{remark}
A result analogous to \ref{8}, combined with a more careful study of the cohomological 
properties of $E_{r,k}$, gives information about effective  
separation of points by the linear series $|k\Theta_{N}|$. We will not insist 
on this aspect in the present article.
\end{remark}

In the same spirit of studying properties of linear series on $U_{X}(r,0)$
via vector bundle techniques, one can look at 
multiplication maps on spaces of sections and normal generation. The 
Verlinde bundles are again an essential tool. The underlying theme is the 
study of surjectivity of the multiplication map
$$H^{0}(\OO(k\Theta_{N}))\otimes H^{0}(\OO(k\Theta_{N}))\overset{\mu_{k}}
{\longrightarrow}H^{0}(\OO(2k\Theta_{N})).$$  
To this respect we have to start by assuming that $k$ is already chosen such 
that $\LL^{k}$ and $E_{r,k}$ are globaly generated (in particular $k\geq r+1$). 
As proved in Theorem \ref{44}, this also 
induces the global generation of $\OO(k\Theta_{N})$. The method will be to 
look at the kernels of various multiplication maps--in the spirit of 
\cite{Lazarsfeld} for example--and study their cohomology vanishing properties. Let 
$M_{k}$ on $U_{X}(r,0)$ and $M_{r,k}$ on $J(X)$ be the vector bundles defined by the sequences:
\begin{equation}
0\longrightarrow M_{k}\longrightarrow H^{0}(\OO(k\Theta_{N}))\otimes \OO
\longrightarrow \OO(k\Theta_{N})\longrightarrow 0
\end{equation}
and 
\begin{equation}
0\longrightarrow M_{r,k}\longrightarrow H^{0}(E_{r,k})\otimes \OO_{J}
\longrightarrow E_{r,k}\longrightarrow 0.
\end{equation}
By twisting (3) with $\OO(k\Theta_{N})$ and taking cohomology, it is clear 
that the surjectivity of $\mu_{k}$ is equivalent to $H^{1}(M_{k}\otimes 
\OO(k\Theta_{N}))=0$. 

On the other hand, the global generation of $\OO(k\Theta_{N})$ implies that 
the natural map ${\rm det}^{*}E_{r,k}\rightarrow \OO(k\Theta_{N})$
is surjective, so we can consider the vector bundle $K$ defined by the 
following sequence:
\begin{equation}
0\longrightarrow K\longrightarrow {\rm det}^{*}E_{r,k}\longrightarrow \OO(k\Theta_{N})\longrightarrow 0. 
\end{equation}
\begin{remark}\label{9}
Fixing $L\in$ Pic$^{0}(X)$, we can also look at the evaluation sequence 
for $\LL^{k}$ on $SU_{X}(r,L)$:
$$0\longrightarrow M_{\LL^{k}}\longrightarrow H^{0}(\LL^{k})\otimes 
\OO_{SU_{X}}\longrightarrow \LL^{k}\longrightarrow 0.$$
The sequence (5) should be interpreted as globalizing this picture. It induces the above sequence 
when restricted to the fiber of the determinant map over each $L$.
\end{remark}
The study of vanishing for $M_{r,k}$ and $K$ will be the key to obtaining 
the required vanishing for $M_{k}$. This is reflected in the top exact sequence 
in the following commutative diagram, obtained from (3),(4) and (5) as an 
application of the snake lemma:
$$\xymatrix{
& & 0 \ar[d] & 0 \ar[d] & &  \\ 
& 0 \ar[r] & {\rm det}^{*}M_{r,k} \ar[r] \ar[d] & M_{k} \ar[r] \ar[d] & 
K \ar[r] & 0 \\
& & {\rm det}^{*}H^{0}(E_{r,k})\otimes \OO \ar[d] \ar[r]^{\cong} & H^{0}(k\Theta_{N})\otimes 
\OO \ar[d] & & & \\
0 \ar[r] & K \ar[r] & {\rm det}^{*}E_{r,k} \ar[r] \ar[d] &  
\OO(k\Theta_{N}) \ar[r] \ar[d] & 0 & \\
& & 0 & 0 & & } $$

\medskip
We again state our result in part (b) of the following theorem in a form that 
allows algorithmic applications. The main ingredient is an effective normal generation 
bound for $E_{r,k}$, which is the content of part (a).

\begin{theorem}\label{46}
(a) The multiplication map 
$$H^{0}(E_{r,k})\otimes H^{0}(E_{r,k})\longrightarrow H^{0}(E_{r,k}^{\otimes 2})$$
is surjective for $k\geq 2r+1$.
\newline
(b) Under the global generation assumptions formulated above, the multiplication 
map $$\mu_{k}:H^{0}(\OO(k\Theta_{N}))\otimes 
H^{0}(\OO(k\Theta_{N}))\longrightarrow H^{0}(\OO(2k\Theta_{N}))$$ 
is surjective 
as long as the multiplication map $H^{0}(\LL^{k})\otimes H^{0}(\LL^{k})
\rightarrow H^{0}(\LL^{2k})$ on $SU_{X}(r)$ is surjective and $k\geq 2r+1$. 
\end{theorem}
\begin{proof}
(a) This is not hard to deal with when $k$ is a multiple of $r$, since we know 
from \ref{12} that $E_{r,k}$ decomposes in a particularly nice way. To tackle 
the general case though, we have to appeal to \ref{41}. Concretely, we have to see
precisely when the skew Pontrjagin product $$E_{r,k}\widehat{*} E_{r,k}\cong
E_{r,k}* (-1_{J})^{*}E_{r,k}$$ is globally generated, and since by the initial choice 
of a theta characteristic $E_{r,k}$ is symmetric, this is the same as the 
global generation of the usual Pontrjagin product $E_{r,k}* E_{r,k}$. 
As an aside, recall from \ref{41}
that this would imply the surjectivity of all the multiplication maps
$$H^{0}(t^{*}_{x}E_{r,k})\otimes H^{0}(E_{r,k})\longrightarrow 
H^{0}(t^{*}_{x}E_{r,k}\otimes E_{r,k})$$
for all $x\in J(X)$. 

The global generation of this Pontrjagin product is in turn another application of 
the general cohomological criterion \ref{40} for vector bundles on abelian varieties. 
We first prove that $E_{r,k}* E_{r,k}$ also has 
a nice decomposition when pulled back by an isogeny, namely this time by 
multiplication by $2r$. Denote by 
$F$ the Fourier transform $\widehat{E_{r,k}}$, so that $\widehat{F}\cong
(-1_{J})^{*}E_{r,k}\cong E_{r,k}$. Then we have the following isomorphisms:
$$E_{r,k}* E_{r,k}\cong \widehat{F}* \widehat{F}\cong 
\widehat{F\otimes F},$$ 
where the second one is obtained by the 
correspondence between the Potrjagin product and the tensor product via 
the Fourier-Mukai transform, as in \ref{30}(2).  
Next, as in the previous sections, we look at the behaviour of our bundle when 
pulled back via certain 
isogenies (cf. \ref{30}(1)):
\begin{equation}
{k_{J}}_{*}(E_{r,k}* E_{r,k})\cong 
\widehat{k_{J}^{*}F\otimes k_{J}^{*}F}.
\end{equation}
In \ref{35} we proved that $$k_{J}^{*}F\cong k_{J}^{*}\widehat{E_{r,k}}
\cong k_{J}^{*}E_{k,r}^{*}
\cong \underset{s_{k,r}}{\bigoplus}\OO_{J}(-kr\Theta_{N}).$$ 
and by plugging this into (6) we obtain
$${k_{J}}_{*}(E_{r,k}* E_{r,k})\cong 
(\underset{s_{k,r}}{\bigoplus}\OO_{J}(-kr\Theta_{N})\bigotimes 
\underset{s_{k,r}}{\bigoplus}\OO_{J}(-kr\Theta_{N}))^{\widehat{}}$$
$$\cong \widehat{\underset{s_{k,r}^{2}}{\bigoplus}\OO_{J}(-2kr\Theta_{N})}.$$
Finally we apply $(2r)_{J}^{*}\circ k_{J}^{*}$ to both sides
of the isomorphism above and use the behavior of the Fourier transform of 
a line bundle when pulled back via the isogeny that it determines (see 
\ref{30}(3)). 
Since $E_{r,k}* E_{r,k}$ is a direct summand in $k_{J}^{*}{k_{J}}_{*}(E_{r,k}* E_{r,k})$,
we obtain the desired decomposition:
\begin{equation}
(2r)_{J}^{*}(E_{r,k}* E_{r,k})\cong \bigoplus\OO_{J}
(2kr\Theta_{N}).
\end{equation}
This allows us to apply a trick analogous to the one used in the proof of 
\ref{20}. Namely (7) implies that
$$(2r)_{J}^{*}(E_{r,k}* E_{r,k}\otimes \OO_{J}(-\Theta_{N}))
\cong \bigoplus\OO_{J}((2kr-4r^{2})\Theta_{N}).$$ Thus if we denote by 
$U_{r,k}$ the vector bundle $E_{r,k}* E_{r,k}\otimes 
\OO_{J}(-\Theta_{N})$ we clearly have:
$$h^{i}(U_{r,k}\otimes \alpha)= 0, ~ \forall\alpha\in{\rm Pic}^{0}(X),~ \forall 
i>0 {\rm~ and~ } \forall k\geq 2r+1.$$ 
Pareschi's criterion \ref{40} immediately gives then 
that $E_{r,k}* E_{r,k}$ is globally generated for $k\geq 2r+1$,
since $$E_{r,k}* E_{r,k}\cong U_{r,k}\otimes \OO_{J}(\Theta_{N}).$$
\newline
(b) We will show the vanishing of $H^{1}(M_{k}\otimes \OO(k\Theta_{N}))$. By the 
top sequence in the diagram preceding the theorem, it is enough to prove that 
$$H^{1}(K\otimes \OO(k\Theta_{N}))=0 ~{\rm ~and~}~ H^{1}({\rm det}^{*}M_{r,k}
\otimes \OO(k\Theta_{N}))=0.$$
First we prove the vanishing of $H^{1}(K\otimes 
\OO(k\Theta_{N}))$. The key point is to identify the pull-back of $K$ by 
the \'etale cover $\tau$ in the diagram 
$$\xymatrix{
SU_{X}(r)\times J(X) \ar[d]_{p_{2}} \ar[r]^{\hspace{7mm} \tau} & U_{X}(r,0) \ar[d]^{{\rm det}} \\
J(X) \ar[r]^{r_{J}} & J(X) } $$
described in section 3. In the pull-back sequence 
$$0\longrightarrow \tau^{*}K\longrightarrow \tau^{*}{\rm det}^{*}E_{r,k}
\longrightarrow \tau^{*}\OO(k\Theta_{N})\longrightarrow 0$$
we can identify $\tau^{*}{\rm det}^{*}E_{r,k}$ with $p_{2}^{*}{r_{J}}^{*}E_{r,k}$
and $\tau^{*}\OO(k\Theta_{N})$ with $\LL^{k}\boxtimes \OO_{J}(kr\Theta_{N})$.
In other words we have the exact sequence
$$0\longrightarrow \tau^{*}K\longrightarrow H^{0}(\LL^{k})\otimes 
\OO_{J}(kr\Theta_{N})\longrightarrow \LL^{k}\boxtimes \OO_{J}(kr\Theta_{N})
\longrightarrow 0,$$
which shows that the following isomorphism holds (cf. \ref{9}):
$$\tau^{*}K\cong M_{\LL^{k}}\boxtimes \OO_{J}(kr\Theta_{N}).$$
Finally we obtain the isomorphism 
$$\tau^{*}(K\otimes \OO(k\Theta_{N}))
\cong (M_{\LL^{k}}\otimes \LL^{k})\boxtimes \OO_{J}(2kr\Theta_{N}).$$
Certainly by the argument mentioned earlier the surjectivity of the 
multiplication map $H^{0}(\LL^{k})\otimes H^{0}(\LL^{k})
\rightarrow H^{0}(\LL^{2k})$ is also equivalent to $H^{1}(M_{\LL^{k}}
\otimes \LL^{k})=0$. The required vanishing is then an easy application 
of the K\"unneth formula.  

The next step is to prove the vanishing of $H^{1}({\rm det}^{*}M_{r,k}
\otimes \OO(k\Theta_{N}))$. From the projection formula we know that 
$$R^{i}{\rm det}_{*}({\rm det}^{*}M_{r,k}\otimes \OO(k\Theta_{N}))\cong M_{r,k}
\otimes R^{i}{\rm det}_{*}\OO(k\Theta_{N})=0 ~{\rm ~for~ all~}i>0$$
since obviously $R^{i}{\rm det}_{*}\OO(k\Theta_{N})=0$ for all $i>0$.
The Leray spectral sequence reduces then our problem to proving the 
vanishing $H^{1}(M_{r,k}\otimes E_{r,k})=0$, which is basically equivalent
to the surjectivity of the multiplication map 
$$H^{0}(E_{r,k})\otimes H^{0}(E_{r,k})\longrightarrow H^{0}(E_{r,k}^{\otimes 2}).$$ 
This is the content of part (a).
\end{proof}

\begin{corollary}
For $k$ as in \ref{46}, $\OO(k\Theta_{N})$ is very ample. 
\end{corollary}
\begin{proof}
Since $\Theta_{N}$ is ample, by a standard argument the assertion is true if the 
multiplication maps 
$$H^{0}(\OO(k\Theta_{N})\otimes H^{0}(\OO(kl\Theta_{N}))\longrightarrow 
H^{0}(\OO(k(l+1)\Theta_{N}))$$
are surjective for $l\geq 1$. For $l=1$ this is proved in the theorem and the case 
$l\geq2$ is similar but easier.
\end{proof}

\begin{remark}
The case of line bundles in the theorem above (i.e. $r=1$) is the statement
for Jacobians of a well known 
theorem of Koizumi (see \cite{Koizumi} and \cite{Sekiguchi}). Applied to
that particular case, the metod of 
proof in Step 2 above is of course implicit in Pareschi's paper \cite{Pareschi}.
\end{remark}

For an effective bound implied by the previous theorem we have to restrict 
ourselves to the case of rank $2$ vector bundles, since to the best of our 
knowledge nothing is known about multiplication maps on $SU_{X}(r)$ for 
$r\geq 3$.

\begin{corollary}\label{61}
For a generic curve $X$ the multiplication map 
$$H^{0}(\OO(k\Theta_{N}))\otimes H^{0}(\OO(k\Theta_{N}))\overset{\mu_{k}}
{\longrightarrow}H^{0}(\OO(2k\Theta_{N}))$$ 
on $U_{X}(2,0)$ is surjective for $k\geq {\rm max}\{5,g-2\}$ and so $\OO(k\Theta_{N})$ 
is very ample for such $k$. 
\end{corollary}
\begin{proof}
This follows by a theorem of Laszlo \cite{Laszlo}, which says that on a 
generic curve the multiplication map $$S^{k}H^{0}(\LL^2)\longrightarrow 
H^{0}(\LL^{2k})$$ on $SU_{X}(2)$ is surjective for $k\geq g-2$. See also 
\cite{Beauville} \S4 for a survey of results in this direction.
\end{proof}

\begin{remark}
A more refined study along these lines gives analogous results in the extended setting 
of higher syzygies and $N_{p}$ properties. We hope to come back to this somewhere else. 
\end{remark}

We would like to end this section with another application to multiplication maps.
Although probably not of the same significance as the previous results, it still 
brings some new insight through the use of methods 
characteristic to abelian varieties. Recall from \ref{47}
that the Picard group of $U_{X}(r,0)$
is generated by $\OO(\Theta_{N})$ and the preimages of line bundles on 
$J(X)$. We want to study ``mixed'' multiplication maps of the form:
\begin{equation}
H^{0}(\OO(k\Theta_{N})\otimes H^{0}({\rm det}^{*}\OO_{J}(m\Theta_{N}))
\overset{\alpha}{\longrightarrow} H^{0}(\OO(k\Theta_{N})\otimes {\rm det}^{*}\OO_{J}
(m\Theta_{N})).
\end{equation}

\begin{proposition}\label{49}  
The multiplication map $\alpha$ in (8) is surjective if $m\geq 2$ and 
$k\geq 2r+1$.  
\end{proposition}
\begin{proof}
By repeated use of the projection formula, from the commutative diagram
$$\xymatrix{
H^{0}(\OO(k\Theta_{N})\otimes H^{0}({\rm det}^{*}\OO_{J}(m\Theta_{N})) 
\ar[d]_{\cong} \ar[r]^{\alpha} 
& H^{0}(\OO(k\Theta_{N})\otimes {\rm det}^{*}\OO_{J}(m\Theta_{N})) \ar[d]^{\cong} \\
H^{0}(E_{r,k})\otimes H^{0}(\OO_{J}(m\Theta_{N})) \ar[r]^{\beta} & 
H^{0}(E_{r,k}\otimes \OO_{J}(m\Theta_{N})) } $$
we see that it is enough to prove the surjectivity of the multiplication 
map $\beta$ on $J(X)$.

This is an application of the cohomological criterion \ref{42} going back to Kempf 
\cite{Kempf}. What we need to check is that
$$h^{i}(E_{r,k}\otimes \OO_{J}(l\Theta_{N})\otimes \alpha)=0,~\forall i>0,~\forall
l\geq -2 {\rm ~and~} \forall \alpha \in {\rm Pic}^{0}(J(X)).$$
It is again enough to prove this after pulling back by multiplication by $r$. Now:
$$r_{J}^{*}(E_{r,k}\otimes \OO_{J}(l\Theta_{N})\otimes \alpha)\cong 
\underset{s_{r,k}}{\bigoplus}\OO_{J}((kr+lr^{2})\Theta_{N})\otimes r_{J}^{*}\alpha$$
hence the required vanishings are obvious as long as $l\geq -2$ and $k\geq 2r+1$. 
\end{proof}

\section{\textbf{Variants for arbitrary degree}}
 
As it is natural to expect, some of the facts discussed in the previous sections for 
moduli spaces of vector bundles of degree $0$ can be extended to arbitrary degree. On the 
other hand, as the reader might have already observed, there are results that do not admit
(at least straightforward) such extensions. In this last paragraph we would like to 
emphasize what can and what cannot be generalized using the present techniques.

Fix $r$ and $d$ arbitrary positive integers. Then we can look at the moduli space $U_{X}(r,d)$
of semistable vector bundles of rank $r$ and degree $d$. For $A\in {\rm Pic}^{d}(X)$, denote 
also by $SU_{X}(r,A)$ the moduli space of rank $r$ bundles with fixed determinant $A$. 
We will write $SU_{X}(r,d)$ when it is not important what determinant is involved.
On these 
moduli spaces one can construct generalized theta divisors as in the degree $0$ case. More 
precisely, denote
$$h={\rm gcd}(r,d),~r_{1}=r/h {\rm ~and~} d_{1}=d/h.$$
Then for any vector bundle $F$ of rank $r_{1}$ and degree $d_{1}$, we can consider $\Theta_{F}$ 
to be the closure in $U_{X}(r,d)$ of the locus 
$$\Theta_{F}^{s}=\{E|~h^{0}(E\otimes F)\neq 0\}\subset U_{X}^{s}(r,d).$$
This does not always have to be a proper subset, but it is so for generic $F$ (see 
\cite{Hirschowitz}) and in that case $\Theta_{F}$ is a divisor. We can of course do the same 
thing with vector bundles of rank $kr_{1}$ and degree $kd_{1}$ for any $k\geq 1$ and a formula
analogous to \ref{32} holds. On $SU_{X}(r,d)$ there are similar divisors $\Theta_{F}$ and 
they all determine the same determinant line bundle. As before this generates ${\rm Pic}
(SU_{X}(r,d))$ and is denoted by $\LL$.

To study the linear series determined by these divisors we define Verlinde type bundles as in 
Section 2. They will now depend on two parameters (again not emphasized by the notation). 
Namely fix $F\in U_{X}(r_{1},r_{1}(g-1)-d_{1})$ and $L\in {\rm Pic}^{d}(X)$. Consider the 
composition 
$$\pi_{L}:U_{X}(r,d)\overset{{\rm det}}{\longrightarrow} {\rm Pic}^{d}(X)
\overset{\otimes L^{-1}}{\longrightarrow} J(X)$$
and define:
$$E_{r,d,k}(=E_{r,d,k}^{F,L}):={\pi_{L}}_{*}\OO(k\Theta_{F}).$$
This is a vector bundle on $J(X)$ of rank $s_{r,d,k}:=h^{0}(SU_{X}(r,d),\LL^{k})$.
There is again a fiber diagram
$$\xymatrix{
SU_{X}(r,A)\times J(X) \ar[d]_{p_{2}} \ar[r]^{\hspace{7mm} \tau} & U_{X}(r,d) 
\ar[d]^{{\pi}_{L}} \\
J(X) \ar[r]^{r_{J}} & J(X) }$$
where $\tau$ is given by tensor product and the top and bottom maps are Galois with 
Galois group $X_{r}$. By \cite{Donagi} \S3 one has the formula
$$\tau^{*}\OO(\Theta_{F})\cong \LL\boxtimes \OO_{J}(krr_{1}\Theta_{N}),$$
where $N\in {\rm Pic}^{g-1}(X)$ is a line bundle such that $N^{\otimes r}\cong L\otimes 
({\rm det}F)^{\otimes h}$. As in \ref{1} we obtain the decomposition:
$$r_{J}^{*}E_{r,d,k}\cong \underset{s_{r,d,k}}{\bigoplus}\OO_{J}(krr_{1}\Theta_{N}).$$
The basic duality setup via Fourier-Mukai transform presented in Section 3 can be extended 
with a little extra care to this general setting. The purpose is to relate linear series on
the complementary moduli spaces $SU_{X}(r,d)$ and $U_{X}(kr_{1},kr_{1}(g-1)-kd_{1})$ (cf.
\cite{Donagi} \S5) and this is realized via a diagram of the form: 
$$\xymatrix{
U_{X}(r,d) \ar[dd]_{{\pi}_{L}} & & U_{X}(kr_{1},kr_{1}(g-1)-kd_{1}) \ar[dd]^{{\pi}_{M}} \\
&  J(X)\times J(X) \ar[dl]_{p_{1}} \ar[dr]^{p_{2}} \\
J(X) & & J(X) }$$ 
where $L\in {\rm Pic}^{d}(X)$, $M\in {\rm Pic}^{kr_{1}(g-1)-kd_{1}}(X)$ and $\pi_{L}$ and 
$\pi_{M}$ are defined as above. One can also choose a vector bundle $G\in U_{X}(r_{1},d_{1})$ and 
consider the Verlinde bundle $E_{kr_{1},kr_{1}(g-1)-kd_{1},h}$ associated to $G$ 
and $M$. As before, there is an obvious tensor product map
$$U_{X}(r,d)\times U_{X}(kr_{1},kr_{1}(g-1)-kd_{1})\overset{\tau}{\longrightarrow}
U_{X}(krr_{1},krr_{1}(g-1)).$$
With the extra (harmless) assumption on our choices that $L\cong ({\rm det}G)^{\otimes h}$ 
and $M\cong ({\rm det}F)^{\otimes k}$ we can show exactly as in \ref{36} that
$$\tau^{*}\OO(\Theta)\cong \OO(k\Theta_{F})\boxtimes \OO(h\Theta_{G}) \otimes (\pi_{L}\times 
\pi_{M})^{*}\mathcal{P},$$
where $\mathcal{P}$ is a normalized Poincar\'e line bundle on $J(X)\times J(X)$. Note that 
this is slightly different from \ref{36} in the sense that we are twisting up to slope $g-1$ 
and on $U_{X}(krr_{1},krr_{1}(g-1))$, $\Theta$ represents the canonical theta divisor. The two
formulations are of course equivalent.    

The unique nonzero section of $\OO(\Theta)$ induces then a nonzero map:
$$SD: E_{r,d,k}^{*}\longrightarrow (E_{kr_{1},kr_{1}(g-1)-kd_{1},h})^{\widehat{}}$$
and the global formulation of the full strange duality conjecture is:

\begin{conjecture}(cf. \cite{Donagi} \S5)
$SD$ is an isomorphism.
\end{conjecture}
On the positive side, the properties of the kernel of this map described in \ref{51} still hold.
On the other hand, the method of proof of Theorem \ref{5} cannot be used for arbitrary degree 
even in the level $1$ situation. The point is that these new Verlinde type bundles may 
always fail to 
be simple. Probably the most suggestive example is the case of $r$ and $d$ coprime, when 
the other extreme is attained for any $k$:
\begin{example}
If ${\rm gcd}(r,d)=1$, then $E_{r,d,k}$ decomposes as a direct sum of line bundles for 
any $k$. More precisely:
$$E_{r,d,k}=\underset{s_{r,d,k}}{\bigoplus}\OO(k\Theta_{N}),~ \forall k\geq 1.$$
This can be seen by imitating the proof of \ref{12}.
\end{example}

On a more modest note, the main result of \cite{Donagi} can be naturally integrated into 
these global arguments on $J(X)$. It is obtained by calculus with Fourier transforms in 
the spirit of Section 3 and we do not repeat the argument here:

\begin{proposition}(\cite{Donagi}, Theorem 1)
$h^{0}(U_{X}(r,d),\OO(k\Theta_{F}))=\frac{k^{g}}{h^{g}}\cdot s_{r,d,k}$.
\end{proposition}

Turning to effective global generation and normal generation $U_{X}(r,d)$, the picture 
described in Section 5 completely extends, with the appropriate modifications, to the 
general case. All the effective bounds turn out to depend on the number $h={\rm gcd}(r,d)$.
The global generation result analogous to \ref{44} is formulated as follows:

\begin{theorem}
$\OO(k\Theta_{F})$ is globally generated on $U_{X}(r,d)$ as long as $k\geq h+1$
and $\LL^{k}$ is globally generated on $SU_{X}(r,d)$. Moreover, $\OO(k\Theta_{F})$ is not globally generated for $k\leq h$.
\end{theorem}
In \cite{Popa2} \S4 it is proved that $\LL^{k}$ is globally generated on $SU_{X}(r,d)$
for $k\geq {\rm max}\{\frac{(r+1)^{2}}{4r} h,$ $\frac{r^{2}}{4s} h\}$,
where $s$ is an invariant of the moduli space that we will not define here, but satisfying $s\geq h$
so that in particular $k\geq \frac{(r+1)^{2}}{4}$ always works.
This implies then:

\begin{corollary}  
$\OO(k\Theta_{F})$ is globally generated on $U_{X}(r,d)$ for 
$k\geq {\rm max}\{\frac{(r+1)^{2}}{4r} h,\frac{r^{2}}{4s} h\}$.
\end{corollary}
This again produces optimal results in the case of 
rank $2$ and rank $3$ vector bundles (see \cite{Popa2} (5.7)):

\begin{corollary}
$\OO(2\Theta_{F})$ is globally generated on $U_{X}(2,1)$ and $U_{X}(3,\pm 1)$.
\end{corollary}

In the same vein, the normal generation result \ref{46} can be generalized to:

\begin{theorem}\label{50}
The multiplication map 
$$\mu_{k}:H^{0}(\OO(k\Theta_{F}))\otimes 
H^{0}(\OO(k\Theta_{F}))\longrightarrow H^{0}(\OO(2k\Theta_{F}))$$ 
on $U_{X}(r,d)$ is surjective 
as long as the multiplication map $H^{0}(\LL^{k})\otimes H^{0}(\LL^{k})
\rightarrow H^{0}(\LL^{2k})$ on $SU_{X}(r,d)$ is surjective and $k\geq 2h+1$. 
For such $k$, $\OO(k\Theta_{F})$ is very ample.
\end{theorem}

\begin{corollary}
For $X$ generic $\mu_{k}$ is surjective on $U_{X}(2,1)$ if $k\geq {\rm max}\{3,\frac{g-2}
{2}\}$.
\end{corollary}
\begin{proof}
This is a consequence of \ref{50} and \cite{Laszlo}, where it is proved that on $SU_{X}(2,1)$ 
$$S^{k}H^{0}(\LL)\longrightarrow H^{0}(\LL^{k})$$
is surjective for $k\geq g-2$.
\end{proof}

It is certainly not hard to formulate further results corresponding to \ref{49}. 
We leave this to the interested reader.

\end{document}